\documentclass{amsart}
\title{Traces arising from regular inclusions}
\author{Danny Crytser,  Gabriel Nagy}
\usepackage{amssymb,amsmath,amsthm,verbatim}
\usepackage{graphicx,color}

\def\crytser{\textcolor{blue}}
\def\crytserb{\textcolor{black}}

\newcommand{\T}{\mathbb{T}}
\newcommand{\Z}{\mathbb{Z}}

\newcommand{\pstar}{\phantom{*}}

\theoremstyle{plain}
\newtheorem{theorem}{Theorem}[section]
\newtheorem{proposition}[theorem]{Proposition}
\newtheorem{lemma}[theorem]{Lemma}
\newtheorem{corollary}[theorem]{Corollary}
\newtheorem{fact}[theorem]{Fact}
\newtheorem*{claim}{\sc Claim}

\theoremstyle{definition}
\newtheorem{definition}[theorem]{Definition}

\newtheorem{example}[theorem]{Example}
\theoremstyle{remark}
\newtheorem{remark}[theorem]{\sc Remark}
\newtheorem*{question}{Question}

\newtheorem*{notation}{Notation}
\newtheorem*{notations}{Notations}
\newtheorem*{mycomment}{Comment}

\usepackage{enumerate}
\usepackage{hyperref} 
\usepackage{graphicx,tikz,amssymb,amsmath,amsthm,amsfonts}
\usetikzlibrary{arrows, trees}
\usepackage[all]{xy}
\usetikzlibrary{decorations.markings}
\usetikzlibrary{arrows}
\usepackage{hyperref}

\begin{document}

\begin{abstract}
We study the problem of extending a state on an abelian $C^*$-subalgebra to a tracial state on the ambient $C^*$-algebra. We propose an approach that is well-suited to the case of regular inclusions, in which there is a large supply of normalizers of the subalgebra. Conditional expectations onto the subalgebra give natural extensions of a state to the ambient $C^*$-algebra; we prove that these extensions are tracial states if and only if certain invariance properties of both the state and conditional expectations are satisfied. In the example of a groupoid $C^*$-algebra, these invariance properties correspond to invariance of associated measures on the unit space under the action of bisections. Using our framework, we are able to completely describe the tracial state space of a Cuntz-Krieger graph algebra. Along the way we introduce certain operations called graph tightenings, which both streamline our description and provides connections to related finiteness questions in graph $C^*$-algebras. Our investigation has close connections with the so-called unique state extension property and its variants.
\end{abstract}

\maketitle

\section*{Introduction}

%
A {\em trace\/} on an complex algebra $A$ is a linear functional $\phi:A\to\mathbb{C}$ satisfying
$\phi(xy)=\phi(yx)$ for all $x,y \in A$. If $A$ is a $C^*$-algebra, and the trace $\phi$ is also a state, it is simply called a {\em tracial state}.
%
In this paper we study tracial states on $C^*$-algebras $A$ by reconstructing them from their restrictions to {\em abelian\/} subalgebras $B\subset A$. The material is organized as follows


In Section 1 our approach focuses on the case when
a conditional expectation $\mathbb{E}: A \to B$ exists and the ``candidate'' tracial state on $A$ is $\phi\circ \mathbb{E}$, where
$\phi \in S(B)$. In other words, we focus on states on $A$ \crytserb{that} factor through $\mathbb{E}$; equivalently,
states that vanish on $\text{ker}\,\mathbb{E}$. 
\crytserb{In orther to characterize such states, we identify a certain {\em invariance\/} condition on $\phi$, coupled with a a suitable {\em normalization\/} condition on $\mathbb{E}$ (both conditions employ {\em normalizers\/} of $B$).} 

Section 2 specializes  our investigation to the case of \'{e}tale groupoid $C^*$-algebras, where the natural abelian $C^*$-algebra to consider is $C_0(G^{(0)})$ -- the $C^*$-algebra of continuous functions that vanish 
at $\infty$ on the unit space $G^{(0)}$. In this framework, the invariance conditions treated in Section 1 become measure theoretical in nature. 

In Section 3 we explore the link between the invariance and normalization conditions from
Section 1 and certain state extension properties. 
\crytserb{When the so-called \emph{extension property} holds, the tracial state space of $A$ can be completely described by its restrictions to $B$.}

The paper concludes with Section 4, where the case of {\em graph $C^*$-algebras\/} is fully investigated, using the results proved in the previous sections. Given some directed graph $E$, our main goal is the complete parametrization of the tracial state space of the associated $C^*$-algebra $C^*(E)$, solely in graph theoretical language.
Earlier work in this direction (\cite{Tomforde1}, \cite{PaskRen1}) identified the notion of {\em graph traces\/} as a major ingredient. In many instances, graph traces are not sufficient for exhausting all tracial states, and our analysis shows exactly what additional structure is necessary: {\em cyclical tags on graph traces}. The usage of cyclical tags alone, although necessary, is still insufficient for describing all tracial
states on $C^*(E)$; however, this deficiency can be 
 fixed using graph operations called {\em tightenings}. 


\section{Invariant states on abelian C*-subalgebras}

Following \cite{Kumjian} \crytserb{and} \cite{Renault3}, given
a $C^*$-algebra inclusion $B \subset A$, an element $n \in A$ is said to \emph{normalize} $B$ if $nBn^* \cup n^* B n  \subset B$. The collection of such normalizers is denoted by $N_A(B)$, or simply $N(B)$ when there is no danger of confusion. Clearly $N(B)$ is closed under products and adjoints, and contains $B$. A $C^*$-inclusion $B\subset A$ is said to be {\em regular}, if $N(B)$ generates $A$ as a $C^*$-algebra. (Equivalently, if the span of $N(B)$ is dense in $A$.)

Most of the $C^*$-algebra inclusions $B\subset A$ we are going to deal with in this paper are 
{\em non-degenerate}, in the sense that $B$ contains an approximate unit for $A$. (Of course, if $A$ is unital, then non-degeneracy of $B$ is equivalent to the fact that $B$ contains the unit of $A$.) \crytserb{Note that, }if 
$B \subset A$ is a non-degenerate $C^*$-subalgebra, then $n^* n$, $nn^* \in B$ for any $n \in N(B)$.

\begin{definition}\label{inv-state-def}
Assume $B\subset A$ is a non-degenerate and let $\phi$ be a state on $B \subset A$. 
\begin{enumerate}
\item Given $n \in N(B)$, we say that $\phi$ is \emph{$n$-invariant} if 
\begin{equation}
\forall\, b\in B:\,\,\,\phi(nbn^*)=\phi(n^*nb).
\label{normalized-phi-def}
\end{equation}
\item Given $N_0 \subset N(B)$, we say that $\phi$ is \emph{$N_0$-invariant} if $\phi$ is $n$-invariant for all $n \in \Sigma$. 
\item Lastly, if $\phi$ is $N(B)$-invariant, then we simply say that $\phi$ is \emph{fully invariant}.
\end{enumerate}
The collection of fully invariant states on $B \subset A$ is denoted by $S^{\operatorname{inv}}(B)$.
\end{definition}

\begin{mycomment}
The restriction $\tau|_B$ of any tracial state $\tau\in T(A)$ is clearly a fully invariant state on $B$, so 
we have an affine $w^*$-continuous map
\begin{equation}
T(A)\ni\tau\longmapsto \tau|_B\in S^{\text{inv}}(B).
\label{T-Sinv}
\end{equation}
This paper \crytserb{aims at} understanding when the map
\eqref{T-Sinv} is either surjective, or injective, or both.
\end{mycomment}

The most important features of normalizers and invariant states are collected in Proposition \ref{semigroup} below.
Both in its proof and elsewhere in the paper, we are going to employ the following well known technical results and notations.

\begin{fact}\label{fxx}
Assume $x$ is an element in some $C^*$-algebra $A$.
\begin{itemize}
\item[(i)] For any function $f\in C\left([0,\infty)\right)$, the elements
$f(xx^*),\,f(x^*x)\in \tilde{A}$, given by continuous functional calculus, satisfy the equality
\begin{equation}
xf(x^*x)=f(xx^*)x.
\label{fxx=}
\end{equation}
\item[(ii)] When specializing to the $k^{\text{th}}$ root functions $f(t)=t^{1/k}$, we also have the equalities
\begin{equation}
\lim_{k\to\infty}(xx^*)^{1/k}x=
\lim_{k\to\infty}x(x^*x)^{1/k}=x.
\end{equation}
\item[(iii)] If we fix a double sequence $(f_k^\ell)_{k,\ell=1}^\infty$ of polynomials in one variable, such that
\begin{equation}
\forall\,k\in\mathbb{N}:\,\,\,\lim_{\ell\to\infty}tf_k^\ell(t)=t^{1/k}\text{\em, uniformly on compact $K\subset [0,\infty)$}
\end{equation}
(this is possible by the Stone-Weierstrass Theorem), then: 
\begin{align}
\lim_{k\to\infty}\lim_{\ell\to\infty}xf^\ell_k(x^*x)x^*x&=\lim_{k\to\infty}\lim_{\ell\to\infty}xx^*xf^\ell_k(x^*x)=x,
\label{fklx1}\\
\lim_{k\to\infty}\lim_{\ell\to\infty}f^\ell_k(xx^*)xx^*x&=\lim_{k\to\infty}\lim_{\ell\to\infty}xx^*f^\ell_l(xx^*)x=x.
\label{fklx2}
\end{align} 
\end{itemize}
\end{fact}

\begin{proposition} \label{semigroup}
Let $B \subset A$ be \crytserb{a} non-degenerate abelian $C^*$-subalgebra of a $C^*$-algebra $A$.
\begin{itemize}
\item[(i)] $\overline{n B} = \overline{ B n}$ for all $n \in N(B)$.
\item[(ii)] All states $\phi\in S(B)$ are $B$-invariant.
\item[(iii)] If $\phi\in S(B)$ is $n$-invariant for some $n\in N(B)$, then $\phi$ is also $n^*$-invariant. 
\item[(iv)] If $\phi\in S(B)$ is both $n_1$-invariant and $n_2$-invariant, for some $n_1,n_2\in N(B)$, then $\phi$ is also $n_1 n_2$-invariant. 
\item[(v)] If $N_0 \subset N(B)$ is a sub-$*$-semigroup,  generated as a $*$-semigroup by some subset
$W\subset N(B)$, and $\phi\in S(B)$ is $W$-invariant, then $\phi$ is $N_0$-invariant. 
\item[(vi)] A state $\phi \in S(B)$ is fully invariant if and only if 
\begin{equation}
\forall\, n\in N(B):\,\,\,\phi(nn^*)=\phi(n^*n).
\label{phi-fully-inv-n}
\end{equation} 
\end{itemize}
\end{proposition}

\begin{proof}
(i) 
It suffices to show that for any $n\in N(B)$ and any $b\in B$, we have
$nb\in \overline{Bn}$ and $bn\in\overline{nB}$.
If we fix $n$ and $b$, then using the $f_k^\ell$'s from Fact \ref{fxx}, combined with the commutativity of $B$, we have
\begin{equation}
nb=
\lim_{k\to\infty}\lim_{\ell\to\infty}f_k^\ell(nn^*)nn^*nb
=
\lim_{k\to\infty}\lim_{\ell\to\infty}f_k^\ell(nn^*)nbn^*n.
\label{nB=Bn}
\end{equation}
Since $n$ normalizes $B$, we know that $nbn^*\in B$, so the elements
$b_k^\ell=f_k^\ell(nn^*)nbn^*$ all belong to $B$, and then
\eqref{nB=Bn}, which now simply states that $nb=\lim_{k\to\infty}\lim_{\ell\to\infty}b_k^\ell n$,
 clearly proves that $nb\in \overline{Bn}$. The fact that $bn\in \overline{nB}$ is proved exactly the same way.

(ii)
This is obvious, since $B$ is abelian.

(iii) Take a sequence $\{b_k\} \subset B$ such that $bn = \lim_k  n b_k$. Then 
\[
\phi(n^* b n) = \lim_k \phi(n^* n b_k) = \lim_k \phi(n b_k n^*) = \phi(bnn^*) = \phi(nn^* b).\]

(iv) Suppose that $b \in B$. Take a sequence $\{c_k\} \subset B$ such that $(n_1^* n_1) n_2 = \lim_k n_2 c_k$. Then 
\begin{align*}
\phi(n_1 n_2 b n_2^* n_1^*) &= \phi(n_1^* n_1 n_2 b n_2^*)  = \lim_k \phi(n_2 c_k b n_2^*)=
\lim_k \phi(n_2^* n_2 c_k b) = \\
&=\phi(n_2^* n_1^* n_1 n_2 b) ,
\end{align*}
so that $\phi$ is $n_1 n_2$-invariant. 

Part (v) follows immediately from (iii) and (iv).

(vi) The ``if'' implication (for which it suffices to prove \eqref{normalized-phi-def} only for positive $b$)
 follows from the observation, that for any $n\in N(B)$ and any $b\in B^+$, the element 
$x=nb^{1/2}$ is again in $N(B)$, so applying condition \eqref{phi-fully-inv-n} to $x$ will clearly imply
$$\phi(nbn^*)=\phi(b^{1/2}n^*nb^{1/2})=\phi(n^*nb).$$

Conversely, if $\phi$ is fully invariant, then
$$\forall\,n\in N(B):\,\,\phi(nn^*)=\lim_\lambda\phi(nu_\lambda n^*)=\lim_\lambda
\phi(n^*n u_\lambda)=\phi(n^*n),$$
where
$(u_\lambda) \subset B$ be an approximate identity for $A$.
\end{proof}
Besides the notion of invariance for states on a $C^*$-subalgebra, we will also use the following two additional variants. 
\begin{definition}
Given a state $\psi \in S(A)$, we say that an element $x \in A$ \emph{centralizes} $\psi$ if $\psi(xa)=\psi(ax)$ for all $a \in A$. It is easy to see that the set \[ Z_\psi = \{x \in A: x \text{ centralizes } \psi \} \] is a $C^*$-subalgebra of $A$. (Obviously, $\psi$ is always tracial when restricted to $Z_\psi$. In particular, $\psi$ is tracial on $A$, if and only if its centralizer $Z_\psi$ contains a set that generates $A$ as a $C^*$-algebra.)
\end{definition}
\begin{definition}
If $B \subset A$ is a $C^*$-subalgebra and $n \in N(B)$, we will say that a map $\Phi: A \to B$ is \emph{normalized by $n$} if $\Phi(nan^*)=n\Phi(a)n^*$ for all $a \in A$.
\end{definition}

\begin{lemma} \label{centralize}
Let $B \subset A$ be a non-degenerate abelian C*-subalgebra with a conditional expectation $\mathbb{E}: A \to B$, which is normalized by some $n \in N(B)$. For a state $\phi \in S(B)$, the following are equivalent:
\begin{itemize}
\item[(i)] $\phi$ is $n$-invariant state on $B$;
\item[(ii)] $\phi \circ \mathbb{E}\in S(A)$ is a state on $A$, which is centralized by $n$.
\end{itemize} 
\end{lemma}

\begin{proof}
The implication $(ii)\Rightarrow (i)$ is pretty obvious, and holds even without the assumption that $\mathbb{E}$ is normalized by $n$. Indeed, if $b\in B$, then
$nbn^*=\mathbb{E}(nbn^*)$ and
$bn^*n=\mathbb{E}(bn^*n)$, so if $\phi\circ \mathbb{E}$ is centralized by $n$, then:
$$
\phi(nbn^*)=(\phi\circ \mathbb{E})\big(n(bn^*)\big)=
(\phi\circ \mathbb{E})\big((bn^*)n\big)=\phi(bn^*n)=\phi(n^*nb).$$
For the proof of $(i)\Rightarrow (ii)$, we fix $a\in A$ and we show that
$\phi\big(\mathbb{E}(an)\big)=\phi\big(\mathbb{E}(na)\big)$.
Fix polynomials $(f_k^\ell)$ as in Fact \ref{fxx}(iii).
Since $\mathbb{E}$ is a conditional expectation, it follows that
\begin{equation}
\mathbb{E}(an)=\lim_{k\to\infty}\lim_{\ell\to\infty}\mathbb{E}\left(anf_k^\ell(n^*n)n^*n\right)=
\lim_{k\to\infty}\lim_{\ell\to\infty}\mathbb{E}\left(anf_k^\ell(n^*n)\right)n^*n.
\end{equation}
By the $n$-invariance of $\phi$, we have
\begin{align}
\phi\big(\mathbb{E}(an)\big)&= \lim_{k\to\infty}\lim_{\ell\to\infty}\phi\left(\mathbb{E}\left(anf_k^\ell(n^*n)\right)n^*n\right)  = \notag \\
&=\lim_{k\to\infty}\lim_{\ell\to\infty}\phi\left(n\mathbb{E}\left(anf_k^\ell(n^*n)\right)n^*\right).
\end{align}
Because $\mathbb{E}$ is normalized by $n$, with the help of \eqref{fxx=} our computation continues as:
\begin{align} 
\phi\big(\mathbb{E}(an)\big)
&=\lim_{k\to\infty}\lim_{\ell\to\infty}\phi\left(\mathbb{E}(nanf_k^\ell\left(n^*n)n^*\right)\right)=\notag\\
&=
\lim_{k\to\infty}\lim_{\ell\to\infty}\phi\left(\mathbb{E}\left(naf_k^\ell(nn^*)nn^*\right)\right).\label{P-phi-normal}
\end{align}
Since $\mathbb{E}$ is a conditional expectation onto an abelian C*-subalgebra, we have:
\begin{align*}
\mathbb{E}(naf_k^\ell(nn^*)nn^*) =\mathbb{E}(na)f_k^\ell(nn^*)nn^*  =\\= f_k^\ell(nn^*)nn^*\mathbb{E}(na)  = \mathbb{E}(f_k^\ell(nn^*)nn^*na),
\end{align*}
so when we return to \eqref{P-phi-normal} and we also use \eqref{fklx2}, we finally get: 
\begin{equation*}
\phi\big(\mathbb{E}(an)\big) =\lim_{k\to\infty}\lim_{\ell\to\infty} \phi \left( \mathbb{E}\left(f_k^\ell(nn^*)nn^*na\right)\right)=
\phi\big(\mathbb{E}(na)\big). \qedhere 
\end{equation*}
\end{proof}

\begin{theorem}\label{phiP-trace-thm}
Let $B \subset A$ be a non-degenerate abelian C*-subalgebra with a conditional expectation $\mathbb{E}: A \to B$, which is normalized by some set $N_0\subset N(B)$. 
For a state $\phi \in S(B)$, the following are equivalent:
\begin{itemize}
\item[(i)] $\phi$ is $N_0$-invariant;
\item[(ii)] $\phi\circ \mathbb{E}$ is centralized by all elements of the $C^*$-subalgebra $C^*(B\cup N_0)\subset A$;
\item[(iii)] the restriction $(\phi \circ \mathbb{E})|_{C^*(B\cup N_0)}$ is a tracial state on $C^*(B\cup N_0)$.
\end{itemize} 
\end{theorem}
\begin{proof}
$(i)\Rightarrow (ii)$. Assume $\phi$ is $N_0$-invariant. By Lemma \ref{centralize}, we clearly have the inclusion $N_0\subset Z_{\phi\circ \mathbb{E}}$, so (using the fact that $Z_{\phi\circ \mathbb{E}}$ is a $C^*$-subalgebra of $A$) in order to prove statement (ii), it suffices to show that $\phi\circ \mathbb{E}$ is also centralized by $B$, which is pretty clear, since $B$ is abelian.

The implication $(ii)\Rightarrow (iii)$ is trivial, since any state becomes tracial when restricted to its centralizer.

$(iii)\Rightarrow (i)$. Assume $(\phi \circ \mathbb{E})|_{C^*(B\cup N_0)}$ is a tracial. In particular, 
$N_0$ centralizes this restriction, so by Lemma \ref{centralize} (applied to $C^*(B\cup N_0)$ in place of $A$), it again follows that $\phi$ is $N_0$-invariant.
\end{proof}

\section{Invariant states in the \'{e}tale groupoid framework}

The invariance \crytserb{conditions from Section 1} can be neatly described in the context of {\em \'{e}tale groupoid $C^*$-algebras}, which we briefly recall here.
A \emph{groupoid} is a set $G$ along with a subset $G^{(2)} \subset G \times G$ of \emph{composable pairs} and two functions: composition $G^{(2)} \ni (\alpha,\beta)\longmapsto \alpha\beta\in G$ 
and an involution $G\ni \gamma\longmapsto \gamma^{-1}\in G$ (the inversion), such that the following hold: 
\begin{itemize}
\item[(i)] $\gamma(\eta \zeta) = (\gamma \eta) \zeta$ whenever $(\gamma,\eta),(\eta,\zeta) \in G^{(2)}$; 
\item[(ii)] $(\gamma,\gamma^{-1}) \in G^{(2)}$ for all $\gamma \in G$, and $\gamma^{-1}(\gamma \eta) = \eta$ and $(\gamma \eta) \eta^{-1} = \gamma$ for $(\gamma,\eta) \in G^{(2)}$.  
\end{itemize}
Elements satisfying $u = u^2 \in G$ are called \emph{units} of $G$ and the set of all such units is denoted $G^{(0)} \subset G$ and called the \emph{unit space} of $G$. There are maps $r,s: G \to G^{(0)}$ defined by 
\[
r(\gamma) = \gamma \gamma^{-1} \qquad \qquad s(\gamma) =\gamma^{-1} \gamma \]
that are called, respectively, the \emph{range} and \emph{source} maps. 
If $A,B \subset G$, then 
$$AB = \{\gamma \in G: \exists \alpha \in A, \beta \in B\text{, such that }\alpha \beta = \gamma\}.$$ 
It is not difficult to show that $(\alpha,\beta) \in G^{(2)}$ if and only if $s(\alpha)=r(\beta)$. For a given unit $u \in G^{(0)}$ there is an associated group $G(u) = \{\gamma \in G: r(\gamma) = s(\gamma) = u \}$; this is called the \emph{isotropy} or \emph{stabilizer group} of $u$. The union of all isotropy groups in $G$ forms a subgroupoid of $G$ called $\operatorname{Iso}(G)$, the \emph{isotropy bundle} of $G$. A groupoid is called \emph{principal} (or an \emph{equivalence relation}) if $\operatorname{Iso}(G) = G^{(0)}$; that is, if no unit has non-trivial stabilizer group. 

%
Throughout this present paper a groupoid $G$ will be called \emph{\'etale}, if it is endowed with a Hausdorff, locally compact and second countable topology so that 
\begin{itemize}
\item[(a)] the composition and inversion operations are continuous (the domain of $\circ$ is equipped with the relative product topology), and furthermore,
\item[(b)]  the range and source maps are local homeomorphisms.
\end{itemize}
By condition (b), for each 
$\gamma\in G$, there exists \crytserb{an} open set $\gamma\in X\subset G$, such that
the maps $s(X)\xleftarrow{\,s|_X\,}X
\xrightarrow{\,r|_X\,}r(X)$ are homeomorphisms onto open sets \crytserb{in $G$}; such an $X$ is called a {\em bisection}.
Note that in the \'etale case, the unit space $G^{(0)}$ is in fact {\em clopen\/} in $G$, and all range and source fibers $r^{-1}(u)$, $s^{-1}(u)$, $u\in G^{(0)}$, are discrete in the relative topology; hence 
compact subsets of $G$ intersect any given range (or source) fiber at most finitely many times. 

In order to define a $C^*$-algebra from an \'etale groupoid $G$, it is necessary to specify a $*$-algebra structure on $C_c(G)$. This is given by 
\begin{align*}
(f\times g)(\gamma) &= \sum_{(\alpha,\beta) \in G^{(2)}: \alpha \beta = \gamma} f(\alpha) g(\beta);\\
f^*(\gamma)&=\overline{f(\gamma^{-1})}.
\end{align*}
(Compactness of supports ensures that the sum involved in the definition of the product gives a well-defined element of $C_c(G)$.)
As $G^{(0)}$ is open in $G$, we have an inclusion $C_c(G^{(0)})\subset C_c(G)$, which turns $C_c(G^{(0)})$ into a 
$*$-subalgebra. However, the $*$-algebra operations on $C_c(G^{(0)})$ inherited from
$C_c(G)$ coincide with the usual (pointwise!) operations:
$h^*=\bar{h}$ and $h\times k=hk$, $\forall\,h,k\in C_c(G^{(0)})$. In fact, something similar can be said concerning 
the left and right
$C_c(G^{(0)})$-module structure of $C_c(G)$: for all $f\in C_c(G)$, $h\in C_c(G^{(0)})$ we have
\begin{align}
(f\times h)(\gamma)&=f(\gamma)h\big(s(\gamma)\big);\\
(h\times f)(\gamma)&=h\big(r(\gamma)\big)f(\gamma).
\end{align} 

Following Renault (\cite{Renault}), for an \'etale groupoid $G$, the full $C^*$-norm on $C_c(G)$ is given as
$$\|f\|=\sup\left\{\big\|\pi(f)\big\|\,:\,\pi\text{ non-degenerate $*$-representation of $C_c(G)$}\right\},$$
and the {\em full groupoid $C^*$-algebra} $C^*(G)$ is defined to be the completion of $C_c(G)$ in the full $C^*$-norm.
When restricted to $C_c(G^{(0)})$, the full $C^*$-norm agrees with the usual sup-norm $\|\cdot\|_\infty$, so by completion,
the embedding $C_c(G^{(0)})\subset C_c(G)$ gives rise to a non-degenerate inclusion
$C_0(G^{(0)})\subset C^*(G)$. At the same time, one can also consider the restriction map, which ends up being a contractive map
$\left(C_c(G),\,\|\cdot\|\right)\ni f\longmapsto f|_{G^{(0)}}\in \left(C_c(G^{(0)}),\,\|\cdot\|_\infty\right)$, so by completion one obtains a
contractive linear map $\mathbb{E}:C^*(G)\to C_0(G^{(0)})$, which is in fact a {\em conditional expectation}. We refer to $\mathbb{E}$ as the {\em natural expectation}. 
Using the KSGNS construction associated with $\mathbb{E}$ (\cite{Lance}) we obtain a
$*$-representation $\pi_{\mathbb{E}}:C^*(G)\to\mathcal{L}\left(L^2\left(C^*(G),\mathbb{E}\right)\right)$, where 
$L^2\left(C^*(G),\mathbb{E}\right)$ is the Hilbert $C_0(G^{(0)})$-module obtained by completing
$C^*(G)$ in the norm given by the inner product
$\langle a|b\rangle_{C_0(G^{(0)})}=\mathbb{E}(a^*b)$. With this representation in mind, the quotient
$C^*(G)/\text{ker}\,\pi_{\mathbb{E}}$ is the so-called {\em reduced\/} groupoid $C^*$-algebra, denoted by $C^*_{\text{red}}(G)$.
An alternative description of the ideal $\text{ker}\,\pi_{\mathbb{E}}$ is to employ the usual GNS-representations
$\pi_{ev_u\circ \mathbb{E}}$, associated with the states $ev_u\circ \mathbb{E}\in S\big(C^*(G)\big)$ that are obtained by composing $\mathbb{E}$
with evaluation maps $ev_u:C_0(G^{(0)})\ni h\longmapsto h(u)\in \mathbb{C}$, $u\in G^{(0)}$. With these (honest) representations in mind, we have
$\text{ker}\,\pi_{\mathbb{E}}=\bigcap_{u\in G^{(0)}}\text{ker}\,\pi_{ev_u\circ \mathbb{E}}$.
As was the case with the full groupoid $C^*$-algebra, after composing with the quotient map
$\pi_{\text{red}}:C^*(G)\to C^*_{\text{red}}(G)$, we still have an embedding $C_c(G)\subset C^*_{\text{red}}(G)$, so we can also view $C^*_{\text{red}}(G)$ as the completion of the convolution $*$-algebra $C_c(G)$ with respect to a (smaller) $C^*$-norm, denoted $\|\cdot\|_{\text{red}}$.
As before, when restricted to $C_c(G^{(0)})$, the norm $\|\,\cdot\,\|_{\text{red}}$ agrees with
$\|\,\cdot\,\|_\infty$, so $C_0(G^{(0)})$ still embeds in $C^*_{\text{red}}(G)$, and furthermore, since the natural expectation 
$\mathbb{E}$ vanishes on $\text{ker}\,\pi_{\mathbb{E}}$, we will have a reduced version of natural expectation, denoted by
 $\mathbb{E}_{\text{red}}:C^*_{\text{red}}(G)\to C_0(G^{(0)})$, which satisfies $\mathbb{E}_{\text{red}}\circ\pi_{\text{red}}=\mathbb{E}$..

As pointed out for instance in \cite{Renault3}, a large supply of normalizers for $C_0(G^{(0)})$ are those elements of the groupoid $C^*$-algebra represented by functions $f \in C_c(G)$ supported
in bisections.
We shall refer to such elements as \emph{elementary normalizers} of $C_0(G^{(0)})$. Note that the collection 
$N_{\text{elem}}\big(C_0(G^{(0)})\big)$ of elementary normalizers, along with $0$, 
is a $*$-subsemigroup of 
$N\left(C_0(G^{(0)})\right)$, and furthermore $N_{\text{elem}}\big(C_0(G^{(0)})\big)$ generate the ambient algebra -- $C^*(G)$ or $C^*_{\text{red}}(G)$ -- as a $C^*$-algebra. Using the embedding of $C_c(G)$ in the groupoid (full or reduced) $C^*$-algebra, we interpret $N_{\text{elem}}\big(C_0(G^{(0)})\big)$ as a subset in $C_c(G)$, namely: 
\begin{equation}
N_{\text{elem}}\big(C_0(G^{(0)})\big)=\bigcup_{X\text{ bisection}}C_c(X)\subset C_c(G).
\label{Nelem-pres}
\end{equation}

\begin{mycomment}
In order to avoid any unnecessary notational complications or duplications, the results and definitions in the remainder of this section are stated only using  the reduced $C^*$-algebra $C^*_{\text{red}}(G)$ as the ambient $C^*$-algebra. However, with only a few explicitly noted exceptions,
by composing with
 the quotient $*$-homomorphism $\pi_{\text{red}}:C^*(G)\to C^*_{\text{red}}(G)$, the same results will hold if we use the full $C^*$-algebra
$C^*(G)$ instead; we leave it to the reader to write down the missing statements corresponding to the full case 
(by simply erasing the subscript ``red'' from the statements). 
\end{mycomment}

The \'etale groupoid framework is particularly convenient because one of the hypotheses in Lemma \ref{centralize} above is automatically satisfied.

\begin{proposition}\label{prop-elem-inv}
The natural conditional expectation $\mathbb{E}_{\text{\rm red}}: C^*_{\text{\rm red}}(G) \to C_0(G^{(0)})$ is normalized by all elementary normalizers. 
In particular, for a state $\phi$ on $C_0(G^{(0)})$, the following are equivalent:
\begin{itemize}
\item[(i)] $\phi$ is an $N_{\text{\rm elem}}\big(C_0(G^{(0)})\big)$-invariant state on $C_0(G^{(0)})$;
\item[(ii)] $\phi\circ \mathbb{E}_{\text{\rm red}}$ is a tracial state on $C^*_{\text{\rm red}}(G)$.
\end{itemize}
\end{proposition}
\begin{proof}
Assume $n \in C_c(X)$, for some bisection $X\subset G$. In order to prove the first assertion,
we must show that $\mathbb{E}_{\text{\rm red}}(n\times f\times n^*)=n\times \mathbb{E}_{\text{\rm red}}(f)\times n^*$, for all $f \in C_c(G)$. Fix $f$, as well as $x \in G^{(0)}$. Then 
\[
\mathbb{E}_{\text{\rm red}}(n\times f\times n^*)(u) = \begin{cases} 
|n(\gamma)|^2 f(s(\gamma)) & \text{if }\exists \gamma \in X \cap r^{-1}(u) \cap s^{-1}(\operatorname{supp }f)  \\
0 & \text{ else}
\end{cases} .\]
It is straightforward to verify that this is the same as $\left(n\times \mathbb{E}_{\text{\rm red}}(f)\times n^*\right)(u)$. 

The second statement is a direct consequence of Theorem \ref{phiP-trace-thm}, combined with the fact that
$N_{\text{elem}}(C_0(G^{(0)}))$ generates $C^*_{\text{red}}(G)$ as a $C^*$-algebra.
\end{proof}

We want to characterize the $N_{\text{elem}}(C_0(G^{(0)}))$-invariant states on $C_0(G^{(0)})$ -- hereafter referred to as {\em elementary invariant\/} states -- completely in measure-theoretical terms on $G$. We introduce the following terminology in parallel with Definition \ref{inv-state-def}.  

\begin{definition}
Let $G$ be an \'{e}tale topological groupoid with unit space $G^{(0)}$, and let $\mu$ be a 
positive Radon measure on $G^{(0)}$. 
\begin{enumerate}
\item Given an open bisection $X\subset G$, we say that $\mu$ is \emph{$X$-balanced} if $\mu(XBX^{-1}) = \mu(s(X) \cap B)$ for any Borel set $B \subset G^{(0)}$. 
\item If $\mathcal{X}$ is a family of open bisections, then we say that $\mu$ is \emph{$\mathcal{X}$-balanced} if $\mu$ is $X$-balanced for all $X \in \mathcal{X}$. 
\item If $\mu$ is $X$-balanced for every open bisection $X$, then we say that $\mu$ is \emph{totally balanced}.
\end{enumerate}
\end{definition}

%

\begin{notations}
Given a proper continuous function between locally compact spaces $h: X \to Y$, and a Radon measure 
$\mu$ on $X$, we denote its $h$-pushforward by $h_* \mu$. This is a Radon measure on $Y$, given by $(h_* \mu)(A) = \mu(h^{-1}(A))$, for any Borel set $A \subset Y$. 
Note that the pushforward construction is covariant: $(g \circ f)_* \mu = g_* (f_* \mu)$. 

By Riesz's Theorem,  we have a bijective correspondence 
\begin{equation}
\text{Prob}(X)\ni \mu\longmapsto \phi_\mu \in S\big(C_0(X)\big)
\label{riesz-cor}
\end{equation}
between the space of {\em \crytserb{Radon probability} measures on $X$} and the {\em state space of $C_0(X)$}, defined as follows. For each 
$\mu\in\text{Prob}(X)$, the associated state $\phi_\mu\in S\big(C_0(X)\big)$ is:
$$\phi_\mu(f)=\int _X f(x)\,d\mu(x),\,\,\,f\in C_0(X).$$
On the level of positive linear functionals, the pushforward construction corresponds to {\em composition\/}:
$$(h_*\phi)(f)=\phi\big(f\circ h),\,\,\,f\in C_0(Y),\,h:X\to Y.$$
\end{notations}


\begin{lemma}\label{Xbal-lemma}
With $G$ as above, let $X\subset G$ be an open bisection. For a finite Radon measure $\mu$ on $G^{(0)}$, the following are equivalent:
\begin{itemize}
\item[(i)] $\mu|_{s(X)} = \big(s \circ (r|_X)^{-1}\big)_* (\mu|_{r(X)})$;
\item[(ii)] $\mu\big(s(B)\big)=\mu\big(r(B)\big)$, for all Borel subsets $B\subset X$;
\item[(iii)] $\mu\big(s(K)\big)=\mu\big(r(K)\big)$, for all compact subsets $K\subset X$;
\item[(iv)] $\mu$ is $X$-balanced.
\end{itemize}
{\rm (In condition (i) we use the restriction notation for measures: if $\mu$ is a finite Radon measure on $G^{(0)}$
-- thought as a function $\mu:\text{Bor}(G^{(0)})\to [0,\infty)$, and
$D\subset G^{(0)}$ is some open subset, then $\mu|_D$ is the Radon measure on $D$ obtained by restricting $\mu$ to 
$\text{Bor}(D)$.)}
\end{lemma}

\begin{proof}
The equivalence $(i)\Leftrightarrow (ii)$ is trivial, because the maps $s(X)\xleftarrow{\,s|_X\,}X
\xrightarrow{\,r|_X\,}r(X)$ are homeomorphisms onto open sets.

The equivalence $(ii)\Leftrightarrow (iv)$ follows from the observation that, for any Borel set $B\subset G^{(0)}$, the set $B'=X\cap s^{-1}(B)\subset X$ is Borel, and furthermore, the sets that appear in the definition of $X$-invariance are precisely
$
XBX^{-1}= r(B') 
$ and
$s(X)\cap B = s(B')$.

Lastly, the equivalence $(ii)\Leftrightarrow (iii)$ follows from regularity and finiteness of $\mu$.
\end{proof}

\crytserb{We are interested in balanced measures because they are tied up with 
elementary invariance.}

\begin{lemma} \label{grpdmeasX}
Let $G$ be an \'{e}tale groupoid with unit space $G^{(0)}$, let $\mu$ be a \crytserb{Radon probability} measure on $G^{(0)}$, and let $\phi_\mu$ be the state on the $C^*$-subalgebra $C_0(G^{(0)})\subset C^*_{\text{\rm red}}(G)$ given by  \eqref{riesz-cor}. For an open bisection $X\subset G$, the following conditions are equivalent:
\begin{itemize}
\item[(i)] $\mu$ is $X$-balanced;
\item[(ii)] $\phi_\mu$ is $C_c(X)$-invariant. (As in \eqref{Nelem-pres}, $C_c(X)\subset 
N_{C^*_{\text{\rm red}}(G)}\left(C_0(G^{(0)})\right)$.)
\end{itemize} 
\end{lemma}

\begin{proof}
The entire argument will be based on the following

\begin{claim}
For any $n\in C_c(X)$ and any $b\in C_c(G^{(0)})$, one has the equalities:
\begin{align}
&\phi_\mu(n^*\!\times\! n\!\times\! b)=\int_{s(X)}\!\left|\left(n\circ (s|_X)^{-1}\right)(u)\right|^2 b(u)\,d\left(\mu|_{s(X)}\right)(u);
\label{grpdmeasXcl1}\\
&\phi_\mu(n\!\times\! b\!\times\! n^*)=\int_{r(X)}\!\left|\left(n\circ (r|_X)^{-1}\right)(u)\right|^2\left(b\circ s\circ (r|_X)^{-1}\right)(u)\,d\left(\mu|_{r(X)}\right)(u);
\label{grpdmeasXcl2}\\
&\phi_\mu(n\!\times\! b\!\times\! n^*)=\int_{s(X)}\!\left|\left(n\circ (s|_X)^{-1}\right)(u)\right|^2 b(u)\,d\left(s\circ (r|_X)^{-1}\right)_*\left(\mu|_{r(X)}\right)(u).\label{grpdmeasXcl3}
\end{align}
\end{claim}
The equality \eqref{grpdmeasXcl1} follows from the definition of the convolution multiplication and $*$-involution, which yields 
\[
(n^*\times  n)(u) = \begin{cases} 
\left|n\left((s|_X)^{-1}(u)\right)\right|^2 & u \in s(X) \\
0 & u \not\in s(X) \end{cases} 
\] 
so we can multiply the functions $n^*n$ and $b$ to obtain: 
\[
(n^*\times  n\times b)(u) = \begin{cases} 
\left|n\left((s|_X)^{-1}(u)\right)\right|^2 b(u) & u \in s(X) \\
0 & u \not\in s(X). \end{cases} 
\] 
Likewise, the equality in \eqref{grpdmeasXcl2} follows from
\[
(n\times b\times n^*)(u) = \begin{cases} 
\left|n\left(( r|_X)^{-1}(u)\right)\right|^2 \cdot b\left(s\left((r|_X)^{-1}(u)\right)\right) & u \in r(X) \\
0 & u \not\in r(X)
\end{cases}
\]
which implies that the support of $n\times b\times n^*$ is contained in $X (\operatorname{supp } b) X^{-1}\subset r(X)$.
Lastly, the equality between the right-hand sides of \eqref{grpdmeasXcl2} and \eqref{grpdmeasXcl3} follows immediately by applying the definition of the pushforward
\begin{equation}
\int_{s(X)}f\,d\left(s\circ (r|_X)^{-1}\right)_*\left(\mu|_{r(X)}\right)=
\int_{r(X)}\left(f\circ s\circ (r|_X)^{-1}\right)\,d\left(\mu|_{r(X)}\right),
\label{grpdmeasXpush1}
\end{equation}
to functions $f\in C_c\big(s(X)\big)$ of the form:
$f(u)=\left|n\circ \left((s|_X)^{-1}\right)(u)\right|^2 b(u)$.

Having proved the Claim, the implication $(i)\Rightarrow (ii)$ follows from Lemma \ref{Xbal-lemma}, which yields:
\begin{equation}
\forall\,n\in C_c(X),\,b\in C_c(G^{(0)}): \,\,\,\phi_\mu(n^*\!\times\!n\!\times\!b)=\phi_\mu(n\!\times\!b\!\times\!n^*).
\label{grpdmeasX(i)->(ii)}
\end{equation}
By density, \eqref{grpdmeasX(i)->(ii)} holds for all $n\in C_c(X)$, $b\in C_0(G^{(0)})$, thus
$\phi_\mu$ is $n$-invariant for all $n\in C_c(X)$.

As for the implication $(ii)\Rightarrow(i)$, all we have to observe is that, if $\phi_\mu$ is $C_c(X)$-invariant, then \eqref{grpdmeasX(i)->(ii)} is valid, which by the identities \eqref{grpdmeasXcl1} and \eqref{grpdmeasXcl3}, simply state that the equality
\begin{equation}
\int_{s(X)}f\,d\left(s\circ (r|_X)^{-1}\right)_*\left(\mu|_{r(X)}\right)=
\int_{s(X)}f\,d\left(\mu|_{s(X)}\right),
\label{grpdmeasXpush3}
\end{equation}
holds for all functions of the form:
\begin{equation}
f(u)=\left|\left(n\circ (s|_X)^{-1}\right)(u)\right|^2 b(u),\,\,n\in C_c(X),\,b\in C_c(G^{(0)}).
\label{grpdmeasXpush2}
\end{equation}
Since (using a partition of unity argument) the functions of the above form 
linearly span all functions in $C_c\left(s(X)\right)$, the equality
\eqref{grpdmeasXpush3} simply states that  
$$\left(s\circ (r|_X)^{-1}\right)_*\left(\mu|_{r(X)}\right)=\mu|_{s(X)},$$ 
so by Lemma \ref{Xbal-lemma},
it follows that $\mu$ is indeed $X$-balanced.
\end{proof}

Combining Proposition \ref{prop-elem-inv} with Lemma \ref{grpdmeasX}, we now reach the following conclusion.
 
\begin{theorem}\label{grpdmeas}
Let $G$ be an \'{e}tale groupoid with unit space $G^{(0)}$, let $\mu$ be a probability Radon measure on $G^{(0)}$, and let $\phi_\mu$ be the state on the $C^*$-subalgebra $C_0(G^{(0)})\subset C^*_{\text{\rm red}}(G)$ given by  \eqref{riesz-cor}. 
The following conditions are equivalent:
\begin{itemize}
\item[(i)] $\mu$ is totally balanced;
\item[(ii)] $\phi_\mu$ is elementary invariant;
\item[(iii)] $\phi_\mu$ is fully invariant;
\item[(iv)] $\phi_\mu\circ \mathbb{E}_{\text{\rm red}}$ is a tracial state on $C^*_{\text{\rm red}}(G)$.
\qed
\end{itemize}
\end{theorem}


In concrete situations, one would like to check condition (i) from the above Theorem in an ``economical'' way.
To be \crytserb{more} precise, assuming that a given measure $\mu\in\text{Prob}(G^{(0)})$ is $\mathcal{X}$-balanced, for some collection of bisections $\mathcal{X}$, we seek a natural subalgebra on which $\phi_\mu\circ \mathbb{E}_{\text{red}}$ is tracial (as in Theorem \ref{phiP-trace-thm}), and furthermore find criteria on $\mathcal{X}$ 
which ensure that our subalgebra is in fact \crytserb{all of} $C^*_{\text{red}}(G)$. Parts of the Lemma below mimic corresponding statements from Proposition \ref{semigroup}. (Each one of statements (i)--(iii) has an implicit statement built-in: the new sets, such as $X'$, $X^{-1}$ and $X_1X_2$ are always bisections.)

\begin{proposition} \label{balance}
Let $G$ be an \'{e}tale groupoid with unit space $G^{(0)}$ and let $\mu$ be a Radon probability measure on $G^{(0)}$, 
\begin{itemize}
\item[(i)] If $\mu$ is $X$-balanced, for some bisection $X$, then 
$\mu$ is $X'$-balanced, for any open subset $X'\subset X$. 
\item[(ii)] If $\mu$ is $X$-balanced, for some bisection $X$, then $\mu$ is $X^{-1}$-balanced. 
\item[(iii)] If $\mu$ is both $X_1$- and $X_2$-balanced, for two bisections $X_1$, $X_2$, then
 $\mu$ is $X_1 X_2$-balanced.
\item[(iv)] Assume $X$ is an open set, written as a union $X=\bigcup_{j\in J}X_j$ of bisections, such that
$s|_X,r|_X:X\to G^{(0)}$ are injective. Then $X$ is a bisection, and if
$\mu$ is $X_j$-balanced for all $j\in J$, then $\mu$ is $X$-balanced.
\end{itemize}
\end{proposition}

\begin{proof}
Statements (i) and (ii) are trivial from Lemma \ref{Xbal-lemma}. 

Before we prove (iii), we need some clarifications. 
First of all, the set $X_1X_2$ is obtained as the image of the open set
$$X_1\circ X_2=\{(\alpha,\beta)\in X_1\times X_2\,:\, s(\alpha)=r(\beta)\}=X_1\times X_2\cap G^{(2)}\subset G^{(2)}.$$
under composition map
$m:G^{(2)}\to G$.
Secondly, 
by the bisection property, the restrictions of the coordinate maps
$X_1\xleftarrow{\,\,p_1\,\,}X_1\times X_2
\xrightarrow{\,\,p_2\,\,}X_2$ give rise to two homeomorphisms 
$p_1(X_1\circ X_2)\xleftarrow{\,\,p_1\,\,}X_1\circ X_2
\xrightarrow{\,\,p_2\,\,}p_2(X_1\circ X_2)$ onto open subsets of $X_1$ and $X_2$ respectively, and furthermore
the compositions $s\circ p_1$ and $r\circ p_2$ agree on $X_1\circ X_2$, and the resulting map, denoted here by
$t:X_1\circ X_2\to \subset G^{(0)}$ is a homeomorphism onto an open subset $D\subset G^{(0)}$.
(This open set is simply $D=t(X_1\circ X_2)=s(X_1)\cap r(X_2)$. By construction,
$X_1X_2=\varnothing\,\Leftrightarrow\,s(X_1)\cap r(X_2)=\varnothing$.) 
Furthermore, again by the bisection property, 
$m|_{X_1\circ X_2}:X_1\circ X_2\to X_1X_2$ is also a homeomorphism onto an open set, so composing its inverse with the coordinate maps, we obtain two homeomorphisms $q_k=p_k\circ (m|_{X_1\circ X_2})^{-1}:X_1X_2\to X_k$, $k=1,2$,
which satisfy $s|_{X_1X_2}=s\circ q_1$ and
$r|_{X_1X_2}=r\circ q_2$.
Using all these three homeomorphisms, the fact that $X_1X_2$ is a bisection is obvious. Not only are the maps
$s(X_1X_2)\xleftarrow{\,s|_{X_1X_2}\,}X_1X_2
\xrightarrow{\,r|_{X_1X_2}\,}r(X_1X_2)$ homeomorphisms, but so is the map
$r\circ q_2=s\circ q_1=t\circ (m|_{X_1\circ X_2})^{-1}:X_1X_2\to D$. 

After all these preparations, statement (iii) follows from the observation that
the $X_1$- and $X_2$-balancing features imply that, for any Borel set $B\subset X_1X_2$ we have
\begin{align*}
\mu\big(s(B)\big)=
\mu\big(s\big(q_2(B)\big)\big)=
\mu\big(r\big(q_2(B)\big)\big)= \\=
\mu\big(s\big(q_1(B)\big)\big)=
\mu\big(r\big(q_1(B)\big)\big)=
\mu\big(r(B)\big),
\end{align*}
so the desired conclusion follows from Lemma \ref{Xbal-lemma}.

(iv). Since we have the equalities $s(X)=\bigcup_{j\in J}s(X_j)$ and
$r(X)=\bigcup_{j\in J}s(X_j)$, it follows that $s(X)$ and $r(X)$ are open.
The fact that both $s(X)\xleftarrow{\,s|_X\,}X\xrightarrow{\,r|_X\,}r(X)$
are \crytserb{homeomorphisms} follows by local compactness.

Finally, to prove that $\mu$ is $X$-balanced, we apply criterion (iii) from Lemma \ref{Xbal-lemma}. Start with some compact set $K\subset X$, and using compactness write it as a finite disjoint union $K=\bigcup_{k=1}^nB_{j_k}$, where $B_{j_k}\subset X_{j_k}$, $k=1,\dots,n$ are Borel sets. Using the fact that $\mu$ is $X_j$-balanced for all $j$, we know that
$\mu\big(s(B_{j_k})\big)=
\mu\big(r(B_{j_k})\big)$, for all $k$, so using that $s$ and $r$ are homeomorphisms, we also have
$s(K)=\bigcup_{k=1}^ns(B_{j_k})$
and
$r(K)=\bigcup_{k=1}^nr(B_{j_k})$ (disjoint unions of Borel sets in $s(X)$ and $r(X)$ respectively), so we have
\begin{align*}
\mu\big(s(K)\big)&=\mu\big(\bigcup_{k=1}^ns(B_{j_k})\big)=\sum_{k=1}^n\mu\big(s(B_{j_k})\big)
=\\
&=\sum_{k=1}^n\mu\big(r(B_{j_k})\big)= \mu\big(\bigcup_{k=1}^nr(B_{j_k})\big)=
\mu\big(r(K)\big).\qedhere
\end{align*}
\end{proof}

Using the above result, combined with Lemma \ref{grpdmeasX}, we immediately obtain the following 
\crytserb{measure-theoretic} groupoid analogue of Theorem \ref{phiP-trace-thm}. 

\begin{theorem}\label{thm-C*Gtrace}
Assume $\mathcal{W}$ is a collection of bisections in the \'{e}tale groupoid $G$, 
and let $\mathcal{X}$ be the inverse semigroup generated by $\mathcal{W}$.
For a measure $\mu\in\text{\rm Prob}(G^{(0)})$, the following are equivalent:
\begin{itemize}
\item[(i)] $\mu$ is $\mathcal{W}$-balanced;
\item[(ii)] $\mu$ is $\mathcal{X}$-balanced;
\item[(iii)]  the state
$\phi_\mu\circ \mathbb{E}_{\text{\rm red}}$ is tracial \crytserb{when restricted to} the subalgebra
$$C^*\bigg(C_0(G^{(0)})\cup\bigcup_{W\in\mathcal{W}}C_c(W)\bigg)=
\overline{\text{\rm span}}\bigg(C_0(G^{(0)})\cup\bigcup_{X\in\mathcal{X}}C_c(X)\bigg).\qed$$
\end{itemize}
\end{theorem}

\begin{remark}
A sufficient condition for a collection $\mathcal{X}$ of bisections of $G$ to satisfy the equality
$$\overline{\text{\rm span}}\bigg(C_0(G^{(0)})\cup\bigcup_{X\in\mathcal{X}}C_c(X)\bigg)=
C^*_{\text{\rm red}}(G)$$
is that $\mathcal{X}$ {\em covers $G\smallsetminus G^{(0)}$}. This follows using a standard partition of \crytserb{unity} argument, which implies the equality $C_c(G)=\text{span}\bigg(C_0(G^{(0)})\cup\bigcup_{X\in\mathcal{X}}C_c(X)\bigg)$.
As a consequence, the desired ``economical'' criterion for traciality of $\phi_\mu\circ \mathbb{E}_{\text{red}}$ is as follows.
\end{remark}

\begin{corollary}\label{tot-inv-cor}
Assume $G$, $\mathcal{W}$ and $\mathcal{X}$ are as in Theorem \ref{thm-C*Gtrace}.
If $\mu\in\text{\rm Prob}(G^{(0)})$ is $\mathcal{W}$-balanced, and 
$\mathcal{X}$ covers $G\smallsetminus G^{(0)}$, then
$\phi_\mu\circ \mathbb{E}_{\text{\rm red}}$ is tracial on $C^*_{\text{\rm red}}(G)$.\qed
\end{corollary}

\section{Tracial states via extension properties}

So far, assuming that an non-degenerate abelian $C^*$-subalgebra $B\subset A$ is the range of a conditional expectation $\mathbb{E}: A \to B$, we have examined certain conditions both for a state $\phi\in S(B)$ and for $\mathbb{E}$, that ensure that $\phi\circ \mathbb{E}$ \crytserb{is a trace}. In the groupoid framework, the natural conditional expectation $\mathbb{E}$
\crytserb{exhibited nice behavior} (elementary invariance), so the focus was solely placed on $\phi$. In this section we
provide another framework, in which again the conditional expectation in question will also be normalized by
all $n\in N(B)$. (As a side issue one should also be concerned with the {\em uniqueness\/} of conditional expectation.)

A natural class of subalgebras to which this analysis can be carried on nicely are
Renault's Cartan subalgebras (\cite{Renault3}; see also the Comment following Corollary \ref{phiP-trace-ext1} below).
As it turns out, very little from the Cartan subalgebra machinery is needed for our purposes: the \emph{almost extension property\/} (\cite{NR2}), which requires that the set
\[ P_1(B\uparrow A) = \{\omega \in \hat{B}: \omega \text{ has a unique extension to a state on }A \}\] 
is weak-$*$ dense in $\hat{B}$ -- the Gelfand spectrum of $B$. 
(A slight strengthening of the above condition will be introduced \crytserb{in} the Comment following Lemma 
\ref{lemma-almost} below.)
 
The utility of the almost extension property is exhibited by Lemma \ref{lemma-almost} below, in preparation of which we need the following simple fact.  

\begin{fact}\label{fact-omega}
Let $\omega$ be a state on $B \subset A$ with extension $\theta \in S(A)$, so that $\theta|_B = \omega$. If
 $x,y \in A$ and satisfy either 
\begin{enumerate}
\item $y^* y \in B$ and $\omega(y^*y)=0$, or 
\item $xx^* \in B$ and $\omega(xx^*)=0$, 
\end{enumerate}
then $\theta(xy)=0$. 

In particular, if $b\in B$ satisfies $0\leq b\leq 1$ and $\omega(b)=1$, then
$$\forall\,a\in A:\,\,\theta(a)=\theta(ab)=\theta(ba)=\theta(bab).$$
\end{fact}

\begin{proof}
Apply the Cauchy-Schwarz inequality for the sesquilinear form: 
$$
\langle a | a' \rangle = \theta(a^* a').
$$
The second statement follows from the first one applied with $y=1-b$.
\end{proof}

\begin{lemma}[{compare to \cite[Lemma 6]{Kumjian}}]\label{lemma-almost}
Let $B \subset A$ be a non-degenerate abelian $C^*$-subalgebra with the almost extension property, and let $\mathbb{E}: A \to B$ be a conditional expectation. Then $\mathbb{E}$ is normalized by all $n \in N(B)$.
\end{lemma}

\begin{mycomment}
As noted in \cite{NR2}, the almost extension property implies that at most one conditional expectation $\mathbb{E}: A \to B$ can exist. In the case such an expectation does exist \crytserb{and the almost extension property holds}, we say that the inclusion $B\subset A$ has the
{\em conditional\/} almost extension property.
\end{mycomment}

\noindent{\em Proof of Lemma \ref{lemma-almost}.}
Fix some normalizer $n\in N(B)$, and let us prove that
\begin{equation}
\mathbb{E}(nan^*)=n\mathbb{E}(a)n^*,
\label{lemma-almost-norm}
\end{equation}
for all $a\in A$. Fix polynomials $(f_k^\ell)$ as in Fact \ref{fxx}(iii), so we have 
\begin{equation}
\mathbb{E}(nan^*)=\lim_{k\to\infty}\lim_{\ell\to\infty}\mathbb{E}(nn^*nf_k^\ell(n^*n)a
f_k^\ell(n^*n)n^*nn^*).
\label{lemma-alomost1}
\end{equation}
Likewise, and using also the fact that $\mathbb{E}$ is a conditional expectation, we also have
\begin{align}
n\mathbb{E}(a)n^*&=\lim_{k\to\infty}\lim_{\ell\to\infty}nn^*nf_k^\ell(n^*n)\mathbb{E}(a)
f_k^\ell(n^*n)n^*nn^*=
\notag\\
&=\lim_{k\to\infty}\lim_{\ell\to\infty}n\mathbb{E}(n^*nf_k^\ell(n^*n)a
f_k^\ell(n^*n)n^*n)n^*,
\label{lemma-alomost2}
\end{align}
Inspecting \eqref{lemma-alomost1} and \eqref{lemma-alomost2}, we now see that it suffices to prove
 \eqref{lemma-almost-norm} for elements of the form $a=n^* a_1 n$; in other words, instead of 
\eqref{lemma-almost-norm}, it suffices to prove
\begin{equation}
\forall\,a\in A:\,\,\,\mathbb{E}(nn^*ann^*)=n\mathbb{E}(n^*an)n^*,
\label{lemma-almost-norm2}
\end{equation}
As both sides of this equation belong to $B$, we only need show that 
\[ \omega(\mathbb{E}(nn^* a nn^*))=\omega(n\mathbb{E}(n^* a n) n^*) \tag{*}\] for all $\omega \in P_1(B\uparrow A)$.

Suppose that $\omega(nn^*)=0$. In this case, we have by Fact \ref{fact-omega} that both sides of (*) are zero. Suppose that $\omega(nn^*) > 0$ and define two states $\psi_\omega$ and $\theta_\omega$ on $A$ by 
\[
\psi_\omega(a) = \frac{(\omega\circ \mathbb{E})(nn^*ann^*)}{\omega(nn^*)^2} \text{ and }
\theta_\omega(a) = \frac{\omega(n\mathbb{E}(n^*an)n^*)}{\omega(nn^*)^2} ,\]
so (*) is equivalent to the equality $\psi_\omega=\theta_\omega$ (of states on $A$).
Note that, if $b \in B$, 
then $\psi_\omega(b) = \theta_\omega(b) = \omega(b)$, so that both states $\psi_\omega$ and $\theta_\omega$ are extensions of
$\omega \in P_1(B\uparrow A)$, so by uniqueness we have $\psi_\omega=\theta_\omega$, and (*) is established. \qed
\medskip

In the context of the conditional almost extension property, Theorem \ref{phiP-trace-thm} has the following consequences.

\begin{corollary}\label{phiP-trace-ext1}
Let $B \subset A$ be a non-degenerate abelian C*-subalgebra with the conditional almost extension property, let
$\mathbb{E}: A \to B$ be its (unique) conditional expectation, 
and let $\phi$ be a state on $B$.
\begin{itemize}
\item[$(a)$] For a subset $N_0\subset N(B)$ the following are equivalent:
\begin{itemize}
\item[(i)] $\phi$ is $N_0$-invariant;
\item[(ii)] $\phi\circ \mathbb{E}$ is centralized by all elements of $C^*(B\cup N_0)\subset A$;
\item[(iii)] the restriction $(\phi \circ \mathbb{E})|_{C^*(B\cup N_0)}$ is a tracial state on $C^*(B\cup N_0)$.
\end{itemize} 
\item[$(b)$] In particular, if $B$ is regular, then $\phi\circ \mathbb{E}$ is a trace on $A$ if and only $\phi$ is fully invariant.\qed
\end{itemize}
{\em (Of course, statement $(b)$ can be slightly relaxed, by requiring that $\phi$ is only $N_0$-invariant for a
subset $N_0\subset N(B)$ which together with $B$ generates $A$ as a $C^*$-algebra.)}
\end{corollary}

\begin{mycomment}
A natural class exhibiting the conditional almost extension property are Cartan subalgebras, as defined by Renault in
\cite{Renault3}. They are regular non-degenerate inclusions $B\subset A$, in which
\begin{itemize}
\item $B$ is {\em maximal abelian\/} (masa) in $A$, and
\item there exists a {\em faithful\/} conditional expectation $\mathbb{E}:A\to B$ (which is necessarily unique). 
\end{itemize}
As pointed out for instance in \cite{BNRSW}, Cartan subalgebras do have the the conditional almost extension property, but there are many examples of regular non-degenerate abelian $C^*$-subalgebra inclusions $B\subset A$ with the conditional almost extension property, which are non-Cartan. In fact, for \'etale groupoids, the equivalent condition to the almost extension property is {\em topological principalness\/}: the set of units $u \in G^{(0)}$ with trivial isotropy $G(u)$ is dense in $G^{(0)}$. For topologically principal groupoids, both inclusions $C_0(G^{(0)})\subset C^*_{\text{red}}(G)$ and
$C_0(G^{(0)})\subset C^*(G)$ have the conditional almost extension property. However, since the (full) conditional expectation $\mathbb{E}:C^*(G)\to C_0(G^{(0)})$ is not faithful in general, $C_0(G^{(0)})$ is generally not Cartan in $C^*(G)$.
On the other hand, since the (reduced) expectation $\mathbb{E}_{\text{red}}:C^*_{\text{red}}(G)\to C_0(G^{(0)})$ is faithful, 
$C_0(G^{(0)})$ is Cartan in $C^*_{\text{red}}(G)$.
\end{mycomment}

Up to this point, we have seen that for regular non-degenerate abelian $C^*$-subalgebras $B\subset A$ with the conditional almost extension property, Corollary \ref{phiP-trace-ext1}(b) provides us with an injective $w^*$-continuous affine map
\begin{equation}
S^{\text{inv}}(B)\ni\phi\longmapsto \phi\circ \mathbb{E}\in T(A),
\label{Sinv-T}
\end{equation}
which is a right inverse of the restriction map
\eqref{T-Sinv}; in particular, it follows that for such inclusions, the map \eqref{T-Sinv} is surjective.

\begin{question}
If $B\subset A$ is a regular non-degenerate abelian $C^*$-subalgebra with the conditional almost extension property, under what additional circumstances is the map \eqref{Sinv-T} also surjective?
(If this is the case, this would imply that the restriction map \eqref{T-Sinv} is in fact an affine $w^*$-homeomorphism.)
\end{question}

As the Example below suggests, even in the case of Cartan inclusions, the map \eqref{Sinv-T} may fail to be 
surjective.

\begin{example}
Let $B = C(\overline{\mathbb{D}}) \subset A = C(\overline{\mathbb{D}}) \rtimes_\alpha \Z = C^*(C(\overline{\mathbb{D}}), u)$, where $\alpha$ is rotation of $\mathbb{D}$ by an irrational multiple of $\pi$ and $u$ is the unitary that implements the automorphism in the crossed product. Then $B$ is a Cartan subalgebra as can be directly verified. The conditional expectation is given on the dense set of Laurent polynomials in $u$ by 
\[
\mathbb{E}( \sum f_n u^n) = f_0 .\] 
(It is obvious that $\mathbb{E}(u^n)=0$ for all $n \neq 0$.)
As $0$ is a fixed point under the rotation $\alpha$, we have that $(\operatorname{ev}_0(\cdot)1, \operatorname{id})$ is a covariant representation of $(C(\overline{\mathbb{D}}),\alpha)$ in $C^*(\mathbb{Z})\cong C(\T)$, thus
it induces a $*$-homomorphism $\rho: A \to C(\T)$. Any state $\psi$ on $C(\T)$ defines a state $\psi \circ \rho$ on $A$, which is clearly tracial since $C(\T)$ is abelian and $\rho$ is a $*$-homomorphism. A tracial state of this form factors through $\mathbb{E}$ if and only if it maps $\{u^n\}_{n \neq 0}$ to $0$, so taking for instance $\psi=ev_z$ to be a point evaluation at $z \in \T$, then clearly $(ev_z\circ\rho)(u)=z \neq 0$, so the trace
$\tau=ev_z\circ \rho\in T(A)$ does not belong to the range of the map \eqref{Sinv-T}.
\end{example}

\begin{remark}
In connection with the above example, the reason that the map $\phi \to \phi \circ \mathbb{E}$ fails to be surjective is the fact that the state $\operatorname{ev}_0$ on $C(\overline{\mathbb{D}})$ does not have a unique extension to a state on $C(\overline{\mathbb{D}}) \rtimes \mathbb{Z}$. Such an obstruction can be avoided if we consider
inclusions with the \emph{(honest) extension property}, which are those non-degenerate abelian $C^*$-subalgebra inclusion $B\subset A$ for which \emph{every} pure state on $B$ has a unique extension to a state on $A$. 
As shown in \cite{KS} and \cite{ABG}, the extension property implies the following: 
\begin{itemize}
\item $B$ is maximal abelian;
\item there exists a unique conditional expectation $\mathbb{E}: A \to B$
\item $\ker \mathbb{E} = [A,B]$ (the closed linear span of the set of elements of the form $ab-ba$, $a \in A, b \in B$). 
\end{itemize}
From the last two properties it follows immediately that any tracial state $\tau \in T(A)$ vanishes on $\ker \mathbb{E}$. Thus, any tracial state factors through $\mathbb{E}$, and is completely determined by its restriction to $B$. 
Since restrictions of the form $\tau\big|_B$, $\tau\in T(A)$ are always fully invariant, Corollary 
\ref{phiP-trace-ext1} has the following immediate consequence. 
\end{remark}

\begin{corollary}\label{cor-Sinv-T-iso}
If $B \subset A$ is a regular abelian $C^*$-subalgebra algebra inclusion with the extension property, 
and $\mathbb{E}: A \to B$ is its associated conditional expectation, then the map
$$S^{\text{\rm inv}}(B)\ni\phi\longmapsto \phi\circ \mathbb{E}\in T(A)$$
is an affine $w^*$-homeomorphism, with inverse  $\tau \to \tau|_B$. \qed
\end{corollary}

%
%


\begin{example} \label{free-group}
For an \'etale groupoid $G$, the inclusions of $C_0(G^{(0)})$ into either the full or reduced $C^*$-algebra of $G$ have the extension property if and only if $G$ is \emph{principal}: all units in $G$ have trivial isotropy group.  In the case when $G$ is a principal groupoid, the above combined with Theorem \ref{grpdmeas} (in both its reduced and full versions) establishes a bijection between the set of totally balanced measures on $G^{(0)}$ and the tracial state spaces of both $C^*(G)$ and $C^*_{\text{red}}(G)$. In particular, if $\Gamma$ is a discrete group acting freely on 
$X$, then the tracial state spaces of both crossed-product $C^*$-algebras $C_0(X) \rtimes \Gamma$ and $C_0(X) \rtimes_{\text{red}} \Gamma$ are naturally identified with the $\Gamma$-invariant Radon probability measures on $X$.  

The condition that the groupoid be principal (or for crossed products, that the action be free) cannot be relaxed, especially in the non-amenable case, as the following example shows. 
Let $\mathbf{F}_2$ -- the free group on two generators --
 act on by translation  on its Alexandrov compactification
$\mathbf{F}_2\cup\{\infty\}$ (by keeping $\infty$ fixed), so that the associated action of
$\mathbf{F}_2$ on the unitized
 on $c_0(\mathbf{F}_2)^\sim$ is given by $\alpha_g(f + c \mathbf{1})= \lambda_g(f) + c \mathbf{1}$, where $\lambda$ is the left-shift action on $c_0(\mathbf{F}_2)$. It is not hard to show that $c_0(\mathbf{F}_2)^\sim \rtimes_{\text{red}} \mathbf{F}_2$ has a unique tracial state. On the other hand, the full crossed product $c_0(\mathbf{F}_2)^\sim \rtimes \mathbf{F}_2$ has the full group $C^*$-algebra $C^*(\mathbf{F}_2)$ as quotient, and so it must have infinitely many tracial states. 

\end{example}

\section{Graph C*-algebras}

In this section we provide a method for parametrizing tracial state spaces on graph $C^*$-algebras. 
\crytserb{Our approach complements the treatment in \cite{Tomforde2} by giving an explicit parametrization of the tracial state space of a graph $C^*$-algebra.}

We begin with a quick review of graph terminology and notation, most of which are borrowed from \cite{Raeburn}.

A \emph{directed graph} $E = (E^0, E^1, r,s)$ consists of two countable sets $E^0,E^1$ as well as range and source maps $r,s: E^1 \to E^0$. A vertex is \emph{regular} if $r^{-1}(v)$ is finite and non-empty. A vertex which is not regular is called \emph{singular}; a singular vertex is either a source ($r^{-1}(v)=\varnothing$) or an infinite receiver ($r^{-1}(v)$ infinite). 

A \emph{finite path} in $E$ is a sequence $\lambda=e_1 \ldots e_n$ of edges satisfying $s(e_k)=r(e_{k+1})$ for $k=1,\ldots,n-1$. (Note that we are using the right-to-left convention.)
The length $\lambda=e_1 \ldots e_n$ is defined to be $|\lambda|=n$, and the set of paths of length $n$ in $E$ is denoted by 
$E^n$; the collection $\bigcup_{n=0}^\infty E^n$ of all finite paths
in $E$ is denoted $E^*$. (The vertices $E^0$ are included in $E^*$ as the paths of length zero.)  
An \emph{infinite path} in $E$ is an infinite sequence $e_1 e_2 \ldots$ of edges in $E$ satisfying $s(e_k)=r(e_{k+1})$ for all $k$; the set of these paths is denoted by $E^\infty$. If $\lambda=e_1 \ldots e_n$ is a finite path then we define its range $r(\lambda)$ to be $r(e_1)$, and its source $s(\lambda)$ to be $s(e_n)$. The range of an infinite path is defined the same way. 
In order to avoid any confusion, for any vertex $v\in E^0$, and any $n\in\mathbb{N}\cup\{\infty\}$,
the set $\{\lambda\in E^n\,:\,r(\lambda)=v,\,|\lambda|=n\}$ will be denoted by $r^{-n}(v)$.

If $\lambda$ is a finite path and $\nu$ is a finite (or infinite) path with $s(\lambda)=r(\nu)$, then we can concatenate the paths to form $\lambda \nu$. Whenever a (finite or infinite) path $\sigma$ can be decomposed as
$\sigma=\lambda\nu$, we write
$\lambda\prec\sigma$ (or $\sigma\succ\lambda$) and we denote $\nu$ by $\sigma\ominus\lambda$. 
A \emph{cycle} is a finite path $\lambda$ of positive length with $r(\lambda)=s(\lambda)$. 

Given a cycle $\lambda=e_1 \ldots e_n\in E^*$,
an \emph{entry} to $\lambda$ is a path $f_1 f_2 \ldots f_j$, $j > 0$, with $r(f_1)=r(e_k)$ and $f_1 \neq e_k$, for some $k$. If no entry to $\lambda$ exists, we say that $\lambda$ is {\em entry-less}. It fairly easy to see that every
entry-less cycle $\lambda$ can be written uniquely as a repeated concatenation
$\lambda=\nu^m$, of a {\em simple\/} entry-less cycle $\nu$, i.e. the number of vertices in $\nu$ equals $|\nu|$. 

\crytserb{An infinite path $x$ is called \emph{periodic} if there exist $\alpha, \lambda \in E^*$, with $s(\alpha)=r(\lambda)=s(\lambda)$, such that $x= \alpha \lambda^\infty$ (that is, $x$ is obtained by following $\alpha$ and then repeating the cycle $\lambda$ forever). If $x = \alpha \lambda^\infty$, and $\lambda$ has minimal length among any cycle in such a decomposition, then the period of $x$ is defined to be $|\lambda|$ and is denoted $\operatorname{per}(x)$. }

\begin{definition}
If $B$ is a $C^*$-algebra then a \emph{Cuntz-Krieger $E$-family} in $B$ is a set $\{S_e, P_v\}_{e \in E^1,v \in E^0}$, where the $S_e$ are partial isometries with mutually orthogonal range projections and the $P_v$ are mutually orthogonal projections which also satisfy: 
\begin{enumerate}[(i)]
\item $S_e^* S_e^{\pstar} = P_{s(e)}$;
\item $S_e^{\pstar} S_e^* \leq P_{r(e)}$; 
\item if $v$ is regular, then $P_v = \sum_{r(e)=v} S_e^{\pstar} S_e^*$.
\end{enumerate} The $C^*$-subalgebra of $B$ generated by $\{S_e,P_v\}_{e \in E^1,v \in E^0}$ is denoted $C^*(S,P)$. 
The graph algebra $C^*(E)$ is the universal $C^*$-algebra generated by a Cuntz-Krieger $E$-family, $C^*(E)=C^*(s,p)$, where 
$\{s_e,p_v\}$ are the \emph{universal generators}. For any Cuntz-Krieger $E$-family $\{S_e,P_v\}_{e \in E^1,v \in E^0}$ there is a unique $*$-homomorphism $\pi_{S,P}: C^*(E) \to C^*(S,P)$ satisfying $\pi_{S,P}(s_e)=S_e$ and $\pi_{S,P}(p_v)=P_v$. 

For an $E$-family $\{S,P\}$ and a finite path $\lambda = e_1 \ldots e_n$ in $E^*$, there is an associated partial isometry $S_\lambda = S_{e_1} S_{e_2} \ldots S_{e_n}$ in $C^*(S,P)$. (If $|\lambda|=0$, so $\lambda$ reduces to a vertex $v\in E^0$, then $S_\lambda=P_v$.) When specializing to $C^*(E)$, we have partial isometries denoted $s^{}_\lambda$, $\lambda\in E^*$. 

By construction, all $s^{}_\lambda\in C^*(E)$, $\lambda\in E^*$ are partial isometries: the source projection of $s^{}_\lambda$ is $s^*_\lambda s^{}_\lambda=p_{s(\lambda)}$; the range projection $s^{}_\lambda s^*_\lambda$ will be denoted from now on by $p^{}_\lambda$.

\end{definition}

As it turns out, \crytserb{one has the equality}
\begin{equation}
C^*(E)=\overline{\operatorname{span}} \{s_\alpha^{\pstar} s_\beta^*: \alpha, \beta	 \in E^*, s(\alpha)=s(\beta) \}.
\label{C*E=span}
\end{equation}
The products $s_\alpha^{\pstar} s_\beta^*$ listed in the right-hand side of \eqref{C*E=span} are referred to as the
\emph{spanning monomials}, and the set of all these elements is denoted by $G(E)$.
The equality \eqref{C*E=span} is due to the fact that $G(E)\cup\{0\}$ is a $*$-semigroup, which is a
consequence of the following product rule:
\begin{equation} 
(s^{\pstar}_\alpha s_\beta^*)(s_\lambda^{\pstar}s_\nu^*)=
\begin{cases}
s^{\pstar}_\alpha s_{\nu(\beta\ominus\lambda)}^*,&\text{if }\lambda\prec\beta\\
s^{\pstar}_{\alpha(\lambda\ominus\beta)}s^*_\nu,&\text{if }\beta\prec\lambda\\
0,&\text{otherwise}
\end{cases}
\label{G-prod}
\end{equation}

Since all projections $p^{}_v$, $v\in E^0$ are mutually orthogonal, for any finite set $V\subset E^0$, the sum
$q^{}_V=\sum_{v\in V}p^{}_v$ will be again a projection, and furthermore, the net
$(q_V^{})_{V\in \mathcal{P}_{\text{fin}}(E^0)}$ forms an approximate unit for $C^*(E)$, hereafter referred to as the 
{\em canonical\/} approximate unit. The $*$-subalgebra
$\bigcup_{V\in\mathcal{P}_{\text{fin}}(E^0)}q^{}_V C^*(E)q^{}_V$ will be denoted by $C^*(E)_{\text{fin}}$.
 
Passing from a graph to a sub-graph does not always produce a meaningful link between the associated $C^*$-algebras.
The best suited objects that allow such links are the identified as follows: given some graph $E$, 
a subset $H \subset E^0$ is called 
\begin{itemize}
\item \emph{hereditary}, if $r(e) \in H$ implies $s(e) \in H$
\item \emph{saturated}, if whenever $v \in E^0$ is regular and $\{s(e): e \in r^{-1}(v) \} \subset H$, it follows 
that $v \in H$.
\end{itemize}
Any subset $H \subset E^0$ is contained in a minimal saturated set $\overline{H}$ called its \emph{saturation}, which
 is the union $\overline{H}=\bigcup_{k=0}^\infty H_k$, where $H_0 = H$ and, for $k > 1$, 
\begin{equation} H_k = H_{k-1}\cup \{v \in E^0:\text{ $v$ regular and } s(r^{-1}(v)) \subset H_{k-1} \}. \label{saturation}\end{equation} 
Clearly, the saturation of a hereditary set is again hereditary. 
The main point about considering such sets is the fact (see \cite{Raeburn}) that, whenever $H\subset E^0$ is saturated and hereditary, and we form the sub-graph
$$E \setminus H = (E^0 \smallsetminus H, s^{-1}(E^0 \smallsetminus H), r,s),$$ 
then we have a natural surjective 
$*$-homomorphism $\rho_H:C^*(E)\to C^*(E\setminus H)$, defined on the generators as
$$
\rho_H(p^{}_v)=\begin{cases}
p^{}_v,&\text{if }v\in E^0\smallsetminus H;\\
0,&\text{otherwise;}\end{cases}
\quad
\rho_H(s^{}_e)=
\begin{cases}
s^{}_e,&\text{if }s(e)\in E^0\smallsetminus H;\\
0,&\text{otherwise.}\end{cases}
$$
(A sub-graph of this form will be called {\em canonical}.)
The ideal $\ker\rho_H$ is simply the closed two-sided ideal generated by $\{p^{}_v\}_{v\in H}$; alternatively, it is also described as:
$$\ker\rho_H=\overline{\text{span}}\{s^{}_\alpha s^*_\beta\,:\,\alpha,\beta\in E^*,\,s(\alpha)=s(\beta)\in H\}.$$

The {\em gauge action\/} on $C^*(E)$ is the point-norm continuous group homomorphism
$\gamma:\mathbb{T}\ni z\longmapsto \gamma_z\in\text{Aut}\big(C^*(E)\big)$, given on the generators by
$\gamma_z(p_v)=p_v$, $v\in E^0$ and $\gamma_z(s_e)=z s_e$, $e\in E^1$. On the spanning monomials listed above, the automorphisms $\gamma_z$, $z\in\mathbb{T}$, act as
$\gamma_z(s_\alpha^{\pstar} s_\beta^*)=z^{|\alpha|-|\beta|}s_\alpha^{\pstar} s_\beta^*$.
The {\em gauge invariant uniqueness theorem\/} of an Huef and Raeburn (see \cite{HuefRaeburn}) states that, 
given some $C^*$-algebra $\mathcal{A}$ equipped with a group homomorphism
$\theta:\mathbb{T}\ni z\longmapsto \theta_z\in\text{Aut}(\mathcal{A})$, and a gauge invariant $*$-homomorphism
$\pi:C^*(E)\to \mathcal{A}$ (that is, such that $\theta_z\left(\pi(x)\right)=
\pi\left(\gamma_z(x)\right)$, $\forall\,x\in C^*(E)$, $z\in\mathbb{T}$), the condition that $\pi$ is injective is equivalent to the condition that $\pi(p_v^{})\neq 0$, for all $v\in E^0$.


There are two distinguished abelian $C^*$-subalgebras of $C^*(E)$ which we use to define states on $C^*(E)$, the first of which is defined as follows.
 
\begin{definition}
Let $E$ be a directed graph. Then the \emph{diagonal} $\mathcal{D} \subset C^*(E)$ is the $C^*$-subalgebra of $C^*(E)$ generated by the set $G_{\mathcal{D}}(E)=\{p^{}_\alpha\}_{\alpha \in E^*}$.
(We sometimes use the notation $\mathcal{D}(E)$, when specifying the graph is necessary.)
\end{definition}

\begin{remark} \label{diag-desc} %
As it turns out, $G_{\mathcal{D}}(E)\cup \{0\}$ is an abelian semigroup of projections; more specifically, by 
\eqref{G-prod}, the product rule for
$G_{\mathcal{D}}(E)$ is:
\begin{equation}
p^{}_\alpha p^{}_\beta=
p^{}_\beta p^{}_\alpha =
\begin{cases}
p^{}_\alpha,&\text{if }\beta\prec\alpha\\
p^{}_\beta,&\text{if }\alpha\prec\beta\\
0,&\text{otherwise}
\end{cases}
\label{GD-prod}
\end{equation}
Using the semigroup property, it follows that we can in fact 
present $\mathcal{D}(E) = \overline{\operatorname{span}}\,G_{\mathcal{D}}(E)$.
We can also write
$\mathcal{D}(E)=\left[\sum_{v\in E^0}\mathcal{D}(E)p^{}_v\right]^{-}$, with each summand presented
as
$$\mathcal{D}(E)p^{}_v=
\overline{\text{span}}\{p^{}_\alpha\,:\,\alpha \in E^*,\,p^{}_\alpha\leq p^{}_v\}
=\overline{\text{span}}\{p^{}_\alpha\,:\,\alpha \in E^*,\,\,r(\alpha)=v\}.$$
As it turns out, each corner $\mathcal{D}(E)p^{}_v$ is in fact a unital abelian AF-subalgebra, with unit $p^{}_v$, so
$\mathcal{D}$ itself is an abelian AF-algebra, which contains the canonical approximate unit 
$(q_V^{})_{V\in \mathcal{P}_{\text{fin}}(E^0)}$.

As explained for instance in \cite{NR}, the Gelfand spectrum $\widehat{\mathcal{D}(E)}$ of the diagonal $C^*$-subalgebra
$\mathcal{D}(E)$ can be identified with the set 
\[ E^{\leq \infty} = E^\infty\cup
\{ x \in E^*: s(x) \text{ is singular } \} \] 
with evaluation maps defined by 
$ev^{\mathcal{D}}_x(p^{}_\alpha) = 1$ if $\alpha\prec x$, and $0$ otherwise. 
In other words, for each $\alpha\in E^*$, when we view $p_\alpha^{}\in\mathcal{D}(E)$ as a continuous function on $\widehat{\mathcal{D}(E)}\simeq E^{\leq\infty}$, this function will be the indicator function of the
compact-open set $Z(\alpha)=\{x\in E^{\leq\infty}\,:\,\alpha\prec x\}$. Furthermore, the sets $Z(\alpha)$, $\alpha\in E^*$ form a basis for the topology, so clearly $\widehat{\mathcal{D}(E)}$ is a totally disconnected.  
When identifying $\mathcal{D}(E)\simeq C_0\big(\widehat{\mathcal{D}(E)}\big)$, the algebraic sum
(without closure!)
$\mathcal{D}(E)_{\text{fin}}=\sum_{v\in E^0}\mathcal{D}(E)p^{}_v$  gets naturally identified with
$C_c\big(\widehat{\mathcal{D}(E)}\big)$, the algebra of continuous functions with compact support. 
\end{remark}

\begin{remark}\label{orthogonality} Cylinder sets can be used to analyze path (in)comparability.
To be more precise, given two paths, $\alpha,\beta\in E^*$, the following statements hold.
\begin{itemize}
\item[I.] (Comparability Rule)
The inequality
$\alpha\prec\beta$ is equivalent to the reverse inclusion $Z(\alpha)\supset Z(\beta)$.
\item[II.] (Orthogonality Rule) Conditions (i)--(iv) below are equivalent:
\begin{itemize}
\item[(i)] $s^*_\alpha s^{}_\beta=0$;
\item[(ii)] the projections $p^{}_\alpha$ and $p^{}_\beta$ are orthogonal, i.e.
$p^{}_\alpha p^{}_\beta=0$;
\item[(iii)] $\alpha$ and $\beta$ are {\em incomparable}, i.e. $\alpha\not\prec\beta$ and
$\beta\not\prec\alpha$;
\item[(iv)] $Z(\alpha)\cap Z(\beta)=\varnothing$.
\end{itemize}
\end{itemize}
\end{remark}

\begin{remark}
Among all paths $x\in E^{\leq\infty}$, the ones of interest to us will be those that represent isolated points in the spectrum
$\widehat{\mathcal{D}(E)}$. 
On the one hand, if $E$ has {\em sources\/} (i.e. vertices $v\in E^0$ with $r^{-1}(v)=\varnothing$), then all 
finite paths that start at sources are determine isolated points in $\widehat{\mathcal{D}(E)}$. 
On the other hand, the infinite paths $x=e_1e_2\dots\in E^\infty$ that produce isolated points in 
$\widehat{\mathcal{D}(E)}$ are precisely those with the property that there exists $k$ such that $r^{-1}(r(e_n))=\{e_n\}$, for all $n\geq k$.
If this is the case, if we form $\alpha=e_1e_2\dots e_{k-1}$, then $\{x\}=Z(\alpha)$. 
Among those paths, the periodic ones will play an important role in our discussion.
\end{remark}

\begin{definition}
A finite path $\alpha=e_1e_2\dots e_n\in E^*$ (possibly of length zero) is called a \emph{ray} if there is a 
a simple entry-less cycle $\nu$, such that $s(\alpha)=s(\nu)$, and furthermore, 
no edge $e_k$ from $\alpha$ is appears in $\nu$. (Note: In \cite{NR}, rays were called {\em distinguished\/} paths.)
In this case, the cycle $\nu$ (which is uniquely determined by $\alpha$) is 
referred to as the {\em seed\/} of $\alpha$. We caution the reader that zero-length rays are permitted: they are what we will call {\em cyclic vertices}. For reasons explained in the second paragraph below, 
the (possibly empty) set of all rays in $E$ will be denoted by $E^*_{\text{\rm\sc ip}}$.

By definition, any two distinct rays $\alpha_1\neq \alpha_2$ are incomparable, so by the Orthogonality Rule (Remark \ref{orthogonality}) they satisfy:
$s^*_{\alpha_1}s^{}_{\alpha_2}=s^*_{\alpha_2}s^{}_{\alpha_1}=0$.

Clearly, rays parametrize the set $E^\infty_{\text{\rm\sc ip}}$ of infinite periodic paths that yield isolated points in
$\widehat{\mathcal{D}(E)}$: 
any such path can be uniquely presented as $x=\alpha\nu^\infty$, with $\alpha$ ray and $\nu$ the seed of $\alpha$, and its period (as a function from $\mathbb{N}$ to $E^1$) is $\text{per}(x)=|\nu|$.
When it would be necessary to emphasize the sole dependence on $\alpha$, we also denote the infinite path 
$\alpha\nu^\infty$ simply by $\xi_\alpha$.
When we collect 
the corresponding points in $\widehat{\mathcal{D}(E)}$, we obtain a countable open set
$\Sigma_{\text{\rm\sc ip}}=\{ev^{\mathcal{D}}_x\,:\,x\in E^\infty_{\text{\rm\sc ip}}\}\subset\widehat{\mathcal{D}(E)}$.
\end{definition}

\begin{remark}\label{path-repn}
Associated with the space $E^{\leq \infty}$ we have the \emph{path representation} $\pi_{\text{path}}: C^*(E) \to B(\ell^2(E^{\leq \infty}))$ given on generators by 
 (see \cite{Raeburn} for details):
$$
\pi_{\text{path}}(s_e)\delta_x = \begin{cases} \delta_{ex} & r(x)=s(e) \\ 0 & \text{ otherwise; } \end{cases} 
\qquad\qquad
\pi_{\text{path}}(p_v)\delta_x  = \begin{cases} \delta_x &  r(x) =v \\ 0 & \text{ otherwise. } \end{cases} 
$$
In general, $\pi_{\text{path}}$ is not faithful; however, it is always faithful on the diagonal subalgebra 
$\mathcal{D}(E)$. This embedding gives us a explicit form of the identification $\widehat{\mathcal{D}(E)} = E^{\leq \infty}$ as follows: for $x \in E^{\leq \infty}$, the associated character on $\mathcal{D}(E)$ is simply $ev^{\mathcal{D}}_x(a)=  \langle \delta_x | \pi_{\text{path}}(a) \delta_x \rangle$. 

For future use, we denote the subalgebras $\pi_{\text{path}}(\mathcal{D}(E))$ and $\pi_{\text{path}}(C^*(E))$ of $B(\ell^2(E^{\leq \infty}))$ by $D_{\text{path}}(E)$ and $A_{\text{path}}(E)$, respectively. 
\end{remark}

%

\begin{notation}
As shown in \cite[Prop. 3.1]{NR}, a spanning monomial $b=s_\alpha^{\pstar} s_\beta^* \in C^*(E)$ is normal if and only if one of the following holds: 
\begin{itemize}
\item[$(a)$] $\alpha = \beta$, so $w=s_\alpha^{\pstar}s_\alpha^*\in G_{\mathcal{D}}(E)$; 
\item[$(b)$] $\alpha\prec\beta$ and $\beta\ominus\alpha$ is an entry-less cycle; 
\item[$(c)$] $\beta\prec\alpha$ and $\alpha\ominus\beta$ is an entry-less cycle. 
\end{itemize}
The set of normal spanning monomials in $C^*(E)$ is denoted by $G_{\mathcal{M}}(E)$.
\end{notation}

\begin{definition}
The \emph{abelian core} $\mathcal{M}(E)$ is the $C^*$-subalgebra of $C^*(E)$ generated by the set $G_{\mathcal{M}}(E)$  of normal spanning monomials.
\end{definition}

\begin{notations}If $b\in G_{\mathcal{M}}(E)\smallsetminus G_{\mathcal{D}}(E)$ 
(i.e. $b$ is of either type $(b)$ or $(c)$ above), 
then $b$ is a normal partial isometry, so its adjoint $b^*$ also acts as its pseudo-inverse. For this reason,
we will denote $b^*$ simply by $b^{-1}$. More generally, we will allow arbitrary negative integer exponents, by letting $b^{-m}$ be an alternative notation for $b^{*m}$. We will also allow zero exponents, by agreeing that 
$b^0=bb^*=b^*b$, a monomial which in fact belongs to $G_{\mathcal{D}}(E)$. (Equivalently, 
for any $b\in G_{\mathcal{M}}(E)\smallsetminus G_{\mathcal{D}}(E)$, the $C^*$-subalgebra $C^*(b)\subset C^*(E)$ generated by $b$ is a unital abelian $C^*$-algebra, and 
$b$ is a unitary element in $C^*(b)$.)
\end{notations}

\begin{remark}
In general, for a monomial $b\in G_{\mathcal{M}}(E)\smallsetminus G_{\mathcal{D}}(E)$, there might be multiple ways to present it as
$s_\alpha^{\pstar} s_\beta^*$, with $\alpha$ and $\beta$ as in $(b)$ or $(c)$ above, but
after careful inspection, one can show that 
$b$ can be uniquely presented as
$b=s_\alpha^{\pstar} s_\nu^m s_\alpha^*=(s_\alpha ^{}s_\nu^{}s_\alpha^*)^m$,
where $\alpha\in E^*$ is a ray with seed $\nu$ and $m$ is some non-zero integer, so if we let $b_\alpha=s_\alpha^{\pstar} s_\nu^{} s_\alpha^*$ (recall that $\nu$ is uniquely determined by $\alpha$), then we can present
$$
G_{\mathcal{M}}(E)\smallsetminus G_{\mathcal{D}}(E)=\{b_\alpha^m\,:\,\text{$\alpha$ ray, $m$ non-zero integer}\}.
$$
Clearly, using our exponent conventions, $G_{\mathcal{M}}(E)\smallsetminus G_{\mathcal{D}}(E)$ is closed under taking adjoints, because $(b_\alpha^m)^*=b_\alpha^{-m}$.
As it turns out, $G_{\mathcal{M}}(E)\cup \{0\}$ is an abelian $*$-semigroup; besides the product rules
\eqref{GD-prod} for $G_{\mathcal{D}}(E)$, the remaining rules which involve the monomials in
$G_{\mathcal{M}}(E)\smallsetminus G_{\mathcal{D}}(E)$ are:
\begin{align}
b_\alpha^0&= p_\alpha^{},\,\,\text{ for all rays }\alpha;\\
b_\alpha^m p_\beta^{}=
p_\beta^{}b_\alpha^m
&=\begin{cases}
b_\alpha^m,&\text{if }\beta\prec\xi_\alpha\\
0,&\text{otherwise}
\end{cases}
\label{GM-prod1}\\
b_{\alpha_1}^{m_1}b_{\alpha_2}^{m_2}=
b_{\alpha_2}^{m_2}b_{\alpha_1}^{m_1}
&=\begin{cases}
b_{\alpha_1}^{m_1+m_2},&\text{if }\alpha_1=\alpha_2\\
0,&\text{otherwise}
\end{cases}
\label{GM-prod2}
\end{align}

By the above $*$-semigroup property, $\mathcal{M}(E)\subset C^*(E)$ is an abelian
$C^*$-subalgebra which contains $\mathcal{D}(E)$, 
and it can also be described as
$\mathcal{M}(E) = \overline{\operatorname{span}}\,G_{\mathcal{M}}(E)$.
Furthermore, the images 
of $\mathcal{D}(E)$ and $\mathcal{M}(E)$ under the path representation agree; that is, $\pi_{\text{path}}(\mathcal{M}(E)) = D_{\text{path}}(E)$. 
In general, $\mathcal{M}(E)$ is much larger than $\mathcal{D}(E)$; in fact,
 $\mathcal{M}(E)=\mathcal{D}(E)'$, the commutant of $\mathcal{D}(E)$ in $C^*(E)$.  

As was the case with the diagonal, we have
$\mathcal{M}(E)=\left[\sum_{v\in E^0}\mathcal{M}(E)p^{}_v\right]^{-}$, with the summand
$\mathcal{M}(E)p^{}_v$ now presented
as
$$\overline{\text{span}}\left(\{b^m_\alpha\,:\,\,m
\in\mathbb{Z},\,\alpha\in E^*_{\text{\rm\sc ip}},\,r(\alpha)=v\}\cup\{p^{}_\alpha\,:\,\alpha\in E^*,\,r(\alpha)=v\}\right),$$
so, upon identifying $\mathcal{M}(E)\simeq C_0\big(\widehat{\mathcal{M}(E)}\big)$, the (non-norm-closed) algebraic sum
$\mathcal{M}(E)_{\text{fin}}=\sum_{v\in E^0}\mathcal{M}(E)p^{}_v$ \crytserb{is} naturally identified with
$C_c\big(\widehat{\mathcal{M}(E)}\big)$, the algebra of continuous functions with compact support.
\end{remark}

\begin{definition}\label{twisted-path-rep}(Twisted path representation.)
With the notation as above, define the \emph{twisted representation} $\Theta: C^*(E) \to 
C\big(\mathbb{T},A_{\text{path}}(E)\big)$ by 
\[
\Theta(a)(z)=\pi_{\text{path}}(\gamma_z(a)) \qquad z \in \mathbb{T}, a \in C^*(E) .\] 
For any pair $(z,x)\in\mathbb{T}\times E^{\leq\infty}$, we define the state
$\omega_{z,x}$ on $C^*(E)$ by
$$\omega_{z,x}(a)=\langle \delta_x | \Theta(a)(z) \delta_x\rangle.$$
\end{definition}

\begin{remark} \label{core-embed}
As $\pi_{\text{path}}$ is injective on $\mathcal{D}(E)$, the gauge-invariant uniqueness theorem implies that $\Theta$ is injective. (The gauge action on the codomain is by translation: $(\lambda_z(f))(w)=f(z^{-1}w)$.) In particular, $\Theta$ yields an injection of $\mathcal{M}(E)$ into $C(\mathbb{T},D_{\text{path}}(E))$. Therefore the spectrum of $\mathcal{M}(E)$ can be recovered as a quotient of the spectrum of $C(\mathbb{T},D_{\text{path}}(E))$ (that is, $\mathbb{T} \times E^{\leq \infty}$), by the natural equivalence relation implemented by $\Theta$. Specifically, if $(z,x) \in \mathbb{T} \times E^{\leq \infty}$, then the restriction
$\omega_{z,x}|_{\mathcal{M}(E)}$ is a pure state on $\mathcal{M}(E)$.
The equivalence relation $\sim$ on 
$\mathbb{T} \times E^{\leq \infty}$ is simply given by:
\begin{equation}
(z_1,x_1)\sim (z_2,x_2)\Leftrightarrow \omega_{z_1,x_1}|_{\mathcal{M}(E)}=\omega_{z_2,x_2}|_{\mathcal{M}(E)}.
\label{equiv-core}
\end{equation}
Since the restrictions of these states on the diagonal act as 
$\omega_{z,x}|_{\mathcal{D}(E)}=ev^{\mathcal{D}}_x$, it is fairly obvious that 
$(z_1,x_1)\sim (z_2,x_2)$ implies $x_1=x_2$.
The precise description of the equivalence classes
$(z,x)_\sim=\{(z_1,x_1)\in \mathbb{T}\times E^{\leq\infty}\,:\,(z_1,x_1)\sim (z,x)\}$ goes as follows.
\begin{equation}
(z,x)_\sim=\begin{cases}z\mathbb{U}_{\text{per}(x)}\times \{x\}_,&\text{if $x\in E^\infty_{\text{\rm\sc ip}}$}\\
\mathbb{T}\times\{x\},&\text{if $x\in E^{\leq\infty}\smallsetminus E^\infty_{\text{\rm\sc ip}}$}
\end{cases}
\end{equation}
(For any integer $n\geq 1$, the symbol $\mathbb{U}_n$ denotes the group of $n^{\text{th}}$ roots of unity.)
\end{remark}


\begin{lemma} \label{top-spec}
Let $E$ be a directed graph.
\begin{itemize}
\item[(i)] When we equip the quotient space $\mathbb{T} \times E^{\leq \infty}\!/\!\sim$ with the quotient topology, the
map $(z,x)_\sim \longmapsto \omega_{z,x}|_{\mathcal{M}(E)}$ is a homeomorphism of  onto the spectrum of $\mathcal{M}(E)$.
\item[(ii)] For every ray $\alpha$, if we regard $p_\alpha^{}$ as a continuous function on $\widehat{\mathcal{M}(E)}$, 
then $p_\alpha^{}$ is the characteristic function of a compact-open subset $\mathbf{T}_\alpha$, which is homeomorphic to $\mathbb{T}$. Specifically, if $\nu$ is the seed of $\alpha$, and $x=\alpha\nu^\infty\in E^\infty_{\text{\rm\sc ip}}$ is the associated periodic path, then
$\mathbf{T}_\alpha=\{(z,x)_\sim\}_{z\in\mathbb{T}}$ and the map
$\mathbb{T}/\mathbb{U}_{|\nu|}\ni z\mathbb{U}_{|\nu|}\longmapsto (z,x)_\sim\in \mathbf{T}_\alpha$ is a homeomorphism.
Alternatively, $\mathbf{T}_\alpha$ is naturally identified with the spectrum -- computed in the unital $C^*$-algebra
$C^*(b_\alpha)$ -- of the normal partial isometry 
$b^{}_\alpha=s_\alpha^{\pstar} s_\nu^{\pstar} s_\alpha^*$. 
\item[(iii)] The compact-open sets 
$(\mathbf{T}_\alpha)_{\alpha\in E^*_{\text{\rm\sc ip}}}$ are mutually disjoint.
When we consider $\Omega_{\text{\rm\sc ip}} = \bigcup_{\alpha\in E^*_{\text{\sc ip}}} \mathbf{T}_\alpha$, and fix a positive Radon 
measure $\mu$ on $\widehat{\mathcal{M}(E)}$ with corresponding positive linear functional $\phi_\mu$ on 
$\mathcal{M}(E)_{\text{\rm fin}}=C_c(\widehat{\mathcal{M}(E)})$, then 
\begin{equation}
\int_{\Omega_{\text{\rm\sc ip}}} f d\mu  = \sum_{\alpha\in E^*_{\text{\rm\sc ip}}} \phi(f p^{}_\alpha) \label{int-cycl}
\end{equation}  for all $f \in \mathcal{M}(E)_{\text{\rm fin}}=C_c(\widehat{\mathcal{M}(E)})$. 
\qed
\end{itemize}
\end{lemma}

\begin{proof}
Parts (i) and (ii) are established in \cite{NR} and \cite{BNR}. For part (iii) we only need to justify the first statement, because the rest follows from the Lebesgue dominated convergence theorem.
This follows immediately from the observation that any two distinct rays $\alpha_1$, $\alpha_2$ are incomparable,
so by \eqref{G-prod} the projections
$p_{\alpha_1}^{}$ and
$p_{\alpha_2}^{}$ are orthogonal, thus the sets 
$\{\mathbf{T}_\alpha\}_{\alpha\text{ ray}}$ form a countable disjoint compact-open cover of $\Omega_{\text{\rm\sc ip}}$. 
\end{proof}

\begin{remark}
Both $\mathcal{D}(E)$ and $\mathcal{M}(E)$ are abelian regular $C^*$-subalgebras in $C^*(E)$, since 
all generators $p_v^{}$, $v\in E^0$ and  $s_e^{}$, $e\in E^1$, normalize both of them.
It is shown in \cite{NR} that $\mathcal{M}(E)$ is in fact a Cartan subalgebra of $C^*(E)$, with its (unique) conditional expectation acting on generators as
\begin{equation}
\mathbb{E}_{\mathcal{M}}(s_\alpha^{\pstar} s_\beta^*) = \begin{cases} 
s_\alpha^{\pstar} s_\beta^*, & \text{if }s_\alpha^{\pstar} s_\beta^* \in G_{\mathcal{M}}(E)\\
0, & \text{otherwise} \\
 \end{cases} \label{pm}
\end{equation}
\end{remark} 

Within this framework, Theorem \ref{phiP-trace-thm} has the following consequence.

\begin{corollary}\label{cor-M-inv}
For a state $\phi$ on $\mathcal{M}(E)$, the following conditions are equivalent:
\begin{itemize}
\item[(i)] The composition $\phi \circ \mathbb{E}_{\mathcal{M}}$ is a tracial state on $C^*(E)$.
\item[(ii)] $\phi$ is $s_e^{}$-invariant for all $e \in E^1$. 
\item[(iii)] $\phi$ is fully invariant.
\qed
\end{itemize}
\end{corollary}

\begin{remark} 
In general, $\mathcal{D}(E)$ is not Cartan, and there may exist more than one conditional expectation onto it.
One expectation -- hereafter referred to as the {\em Haar expectation\/} -- always exists, defined as
$$\mathbb{E}_{\mathcal{D}}(a)=\int_\mathbb{T}\gamma_z\left(\mathbb{E}_{\mathcal{M}}(a)\right)\,dm(z)=
\int_\mathbb{T}\mathbb{E}_{\mathcal{M}}\left(\gamma_z(a)\right)\,dm(z).$$
(Here $m$ denotes the \crytserb{normalized} Lebesgue measure on $\mathbb{T}$; the second equality follows from \eqref{pm}, which clearly implies that
$\mathbb{E}_{\mathcal{M}}$ is gauge invariant.) The Haar expectation acts on the spanning monomials as:
\begin{equation}
\mathbb{E}_{\mathcal{D}}(s_\alpha^{\pstar} s_\beta^*) = \begin{cases} 
p_\alpha^{}, & \text{if }\alpha=\beta\\
0, & \text{otherwise} \\
 \end{cases} \label{pD}
\end{equation}
Since the integration map
$\int_{\mathbb{T}}\gamma_z(a)\,dm(z)$ is always a faithful positive map,
it follows that $\mathbb{E}_{\mathcal{D}}$ is faithful.

Using formulas \eqref{pD} it is easy to see that 
$\mathbb{E}_{\mathcal{D}}$ is also normalized by all $p^{}_v$, $v\in E^0$, and $s^{}_e$, $s^*_e$, 
$e\in E^1$, so we also
have the following analogue of Corollary \ref{cor-M-inv}.
\end{remark} 

\begin{corollary}\label{cor-D-inv}
For a state $\psi$ on $\mathcal{D}(E)$, the following conditions are equivalent:
\begin{itemize}
\item[(i)] The composition $\psi \circ \mathbb{E}_{\mathcal{D}}$ is a tracial state on $C^*(E)$.
\item[(ii)] $\phi$ is $s_e^{}$-invariant for all $e \in E^1$. 
\item[(iii)] $\phi$ is fully invariant.\qed
\end{itemize}
\end{corollary}

\begin{remark}\label{rem-inv}
Either using Corollary \ref{cor-D-inv} or directly \crytserb{from} the definition, it follows that any fully invariant state
$\psi$ on $\mathcal{D}(E)$ satisfies
\begin{equation}
\forall\,\alpha\in E^*:\,\,\,\psi(p^{}_{\alpha})=\psi(p^{}_{s(\alpha)}).
\label{SinvD-proj}
\end{equation}
In particular, a fully invariant state on $\mathcal{D}(E)$ is completely determined by its values on the projections $p_v$, $v\in E^0$.  
\end{remark}

\begin{definition}\label{def-gr-tr}
Let $E$ be a directed graph. A \emph{graph trace} on $E$ is a function $g: E^0 \to [0,\infty)$ such that:
\begin{itemize}
\item[{\sc (a)}] for any $v \in E^0$, $g(v) \geq \sum_{e: r(e)=v} g(s(e))$;
\item[{\sc (b)}] for any regular $v$, we have equality in {\sc (a)}.
\end{itemize}
\crytserb{Note that}, for any graph trace $g$, its null space $N_g=\{v\in E^0\,:\,g(v)=0\}$ is a saturated hereditary set.

Depending on the quantity $\|g\|_1=\sum_{v\in E^0}g(v)$, a graph trace $g$ is declared {\em finite}, 
if $\|g\|_1<\infty$, or 
{\em infinite}, otherwise. 

We denote the set of all graph traces on $E$ by $T(E)$, and the set of finite graph traces on $E$ by
$T_{\text{fin}}(E)$. Lastly, we define the set $T_1(E)=\{g\in T(E):\,\|g\|_1=1\}$, the elements of which are termed
{\em normalized\/} graph traces.
\end{definition}

\begin{theorem}\label{gr-tr-char}
A map $g:E^0\to [0,\infty)$ is a graph trace on $E$, if and only if, every finite 
tuple $\Xi=(\xi_i,\lambda_i)_{i\in I}^n\subset\mathbb{R}\times E^*$ satisfies
\begin{equation}
\textstyle{\sum_{i\in I}\xi_i p^{}_{\lambda_i}\geq 0}
\,\Rightarrow \,
\textstyle{\sum_{i\in I}\xi_ig(s(\lambda_i))\geq 0}.
\label{gr-tr-char-thm}
\end{equation}
\end{theorem}
\begin{proof}
To prove the ``if'' implication, assume $g$ satisfies condition \eqref{gr-tr-char-thm} and let us verify conditions {\sc (a)} 
and {\sc (b)} from Definition \ref{def-gr-tr}. To check condition {\sc(a)}, start off by fixing some $v\in E^0$, and notice that, since for every finite
set $F\subset r^{-1}(v)$, we have $p^{}_v\geq\sum_{e\in F}p^{}_e$ (by the Cuntz-Krieger relations), then by 
\eqref{gr-tr-char-thm}, it follows that $g(v)\geq\sum_{e\in F}g(s(e))$; this clearly implies the inequality
$g(v)\geq\sum_{e\in r^{-1}(v)}g(s(e))$. In order to check {\sc (b)}, simply notice that, if $v$ is regular (so $r^{-1}(v)$ is both finite and non-empty), the by the Cuntz-Krieger relations, we have an equality
$p^{}_v = \sum_{e\in r^{-1}(v)}p^{}_e$, so applying \eqref{gr-tr-char-thm} both ways 
(writing the equality as two inequalities), we clearly get
$g(v)=\sum_{e\in r^{-1}(v)}g(s(e))$.

To prove the ``only if '' implication, we fix a graph trace $g$ and we 
prove the implication \eqref{gr-tr-char-thm}.
As a matter of terminology, if a tuple $\Xi$ satisfies the inequality
\begin{equation}
\textstyle{\sum_{i\in I}\xi_ip^{}_{\lambda_i}\geq 0},
\label{lemma-ineq-tr}
\end{equation}
we will call $\Xi $ {\em admissible}.
Our proof will use induction on the number $\langle \Xi \rangle=|I|+\sum_{i\in I}|\lambda_i|$.

If $\langle \Xi\rangle =1$, then $|I|=1$, thus $I$ is a singleton $\{i_0\}$ and $\lambda_{i_0}$ is a path of length $0$, i.e. a vertex
$v\in E^0$; in this case, \eqref{gr-tr-char-thm} is same as the implication ``$\xi p^{}_v\geq 0\Rightarrow \xi g(v)\geq 0$,'' which is trivial,
since $g$ takes non-negative values. 

Assume \eqref{gr-tr-char-thm} holds whenever $\langle \Xi\rangle <N$, for some $N>1$, and show 
that \eqref{gr-tr-char-thm} holds when $\langle\Xi\rangle =N$. Fix an admissible tuple
$\xi$ with $\langle\Xi\rangle =N$ (so \eqref{lemma-ineq-tr} is satisfied), 
and let us prove the inequality
\begin{equation}
\textstyle{\sum_{i\in I}\xi_jg(s(\lambda_i))\geq 0},
\label{lemma-ineq-tr-conc}
\end{equation}
If we consider the set $W=\{r(\lambda_i)\,:\,i\in I\}$, then we can split (disjointly)
$I=\bigcup_{v\in W}I_v$, where $I_v=\{i\,:\,r(\lambda_i)=v\}$
and we will have
$$\textstyle{\sum_{i\in I}\xi_j g(s(\lambda_i))=
\sum_{v\in W}\sum_{i\in I_v}\xi_i g(s(\lambda_i))},$$
with each tuple
$\Xi_v=(\xi_i,\lambda_i)_{i\in I_v}$ admissible. 
(This is obtained by multiplying the inequality \eqref{lemma-ineq-tr} by $p^{}_v$.) In the case when $W$ has at least two vertices, we have
$\langle \Xi_v\rangle< \langle \Xi\rangle$, $\forall\,v\in W$, so the inductive hypothesis can be used, and the desired conclusion follows.

Based on the above argument, for the remainder of the proof we can assume that $W$ is a singleton, so we have a vertex 
$v\in E^0$, such that $r(\lambda_i)=v$, $\forall\,i\in I$.
Split $I=I^0\cup I^+$, where
$I^0=\{i\in I\,:\,|\lambda_i|=0\}$ and
$I^+=\{i\in I\,:\,|\lambda_i|>0\}$. Since $W$ is a singleton, the set $I^0$ consists of all $I$ for which $\lambda_i=v$.
The case when $I^+=\varnothing$ is trivial, because that would mean that all $\lambda_i$ will be equal to $v$, so for the remainder of the proof we are going to assume that $I^+\neq\varnothing$. 
With this set-up the hypothesis \eqref{lemma-ineq-tr} reads
\begin{equation}
\textstyle{\big(\sum_{i\in I^0}\xi_i\big)p_v+\sum_{i\in I^+}\xi_ip^{}_{\lambda_i}\geq 0},
\label{lemma-ineq-tr0}
\end{equation}
and the desired conclusion \eqref{lemma-ineq-tr-conc} reads:
\begin{equation}
\textstyle{\big(\sum_{i\in I^0}\xi_i\big)g(v)+\sum_{i\in I^+}\xi_i g(s(\lambda_i))\geq 0}.
\label{lemma-ineq-tr-conc0}
\end{equation}
(In the case when $I^0=\varnothing$, we let $\sum_{i\in I^0}\xi_i=0$.)

Since $I^+$ is non-empty (and finite), we can find a finite non-empty 
set $F\subset E^1$ which allows us to split
$I^+$ as a disjoint union of non-empty sets $I^+=\bigcup_{e\in F}I_e$, where $I_e=\{i\in I\,:\,\lambda_i\succ e\}$.
Using the Cuntz-Krieger relations, it follows that the element $q=\sum_{e\in E}s^{}_es^*_e\in\mathcal{D}$ is a 
projection satisfying $q\leq p^{}_v$, so the difference $q'=p^{}_v-q$ is also a (possibly zero) projection.
In either case, it follows that $q's^{}_{\lambda_i}s^*_{\lambda_i}=0$, $\forall\,i\in I^+$, so when we 
multiply \eqref{lemma-ineq-tr0} by $q'$ we obtain:
\begin{equation}
\textstyle{\big(\sum_{i\in I^0}\xi_i\big) q'\geq 0}.
\label{J0-q'}
\end{equation}
Likewise multiplying \eqref{lemma-ineq-tr0} by each $s^{}_{e}s^*_{e}$ we obtain
$$\textstyle{\big(\sum_{i\in I^0}\xi_i\big) s^{}_{e}s^*_{e}+
\sum_{i\in I_e}\xi_js^{}_{\lambda_i} s^*_{\lambda_i}\geq 0},$$
so if we multiply on the left by $s^*_{e}$ and on the right by $s^{}_{e}$, we obtain:
\begin{equation}
\textstyle{\big(\sum_{i\in I^0}\xi_i\big) p^{}_{s(e)}
+\sum_{i\in I_e}\xi_is^{}_{\lambda_i\ominus e} s^*_{\lambda_i\ominus e}\geq 0}.
\label{Xie-ok}
\end{equation}
For each $e\in F$, we can form the tuple $\tilde\Xi_e=(\xi_i,\tilde{\lambda}_i)_{i\in I^0\cup I_e}$ by letting
$$\tilde{\lambda}_i=\begin{cases} s(e),&\text{if }i\in I^0\\
\lambda_i\ominus e,&\text{if }j\in I_e
\end{cases}
$$
and then \eqref{Xie-ok} shows that all $\tilde{\Xi}_e$ are admissible. Since we obviously have
$\langle \tilde{\Xi}_e\rangle <\langle\Xi\rangle$, by the inductive hypothesis we obtain
$\big(\sum_{i\in I^0}\xi_j\big) g(s(e))
+\sum_{i\in I_e}\xi_i g(s({\lambda_i\ominus e}))\geq 0$, which combined with the obvious equality
$s(\lambda_i\ominus e)=s(\lambda_i)$ yields:
\begin{equation}
\textstyle{\big(\sum_{i\in I^0}\xi_j\big) g(s(e))
+\sum_{i\in I_e}\xi_j g(s(\lambda_i))\geq 0}.
\label{Xie-ok-ind}
\end{equation}
We we sum all these inequalities (over $e\in E$), we obtain:
\begin{equation}
\textstyle{\big(\sum_{i\in I^0}\xi_i\big) \big(\sum_{e\in F}^{}g(s(e))\big)
+\sum_{i\in I^+}\xi_i g(s(\lambda_i))\geq 0}.
\label{Xie-ok-ind-almost}
\end{equation}
Comparing this inequality with the desired conclusion \eqref{lemma-ineq-tr-conc0}, we see that it suffices to show that
\begin{equation}
\textstyle{\big(\sum_{i\in I^0}\xi_i\big)g(v)
\geq
\big(\sum_{i\in I^0}\xi_i\big) \big(\sum_{e\in F}^{}g(s(e))\big)}.
\label{lemma-ineq-tr-conc1}
\end{equation}
The case when $I^0=\varnothing$ is trivial, since both sides will equal zero, so for the remainder,
we can assume
$I^0\neq \varnothing$. 
In the case when $q'=0$, that is, when $p^{}_v=\sum_{e\in F}s^{}_es^*_e$, it follows that $v$ is regular and $F=r^{-1}(v)$, so by condition (ii) in the graph trace definition, it follows that
$g(v)=\sum_{e\in F}g(s(e))$ and again \eqref{lemma-ineq-tr-conc1} becomes an equality. Lastly, in the case when
$q'\neq 0$, we use condition (i) in the graph trace definition, which yields
$g(v)\geq\sum_{e\in F}g(s(e))$; this means that desired inequality would follow once we prove that 
$\sum_{i\in I^0}\xi_i\geq 0$, an inequality which is now (under the assumption that $q'$ is a non-zero projection)
a consequence of \eqref{J0-q'}.
\end{proof}

In preparation Proposition~\ref{tr-infinite} below, which contains two easy applications of Theorem~\ref{gr-tr-char}, we introduce the following terminology.

\begin{definition}
A vertex $v\in E^0$ is said to be {\em essentially left infinite}, if there exists an infinite set
$X\subset E^*$ of mutually incomparable paths such that $s(\alpha)=v$ for all $\alpha \in X$. 
\end{definition}

\begin{remark}\label{rem-entry-cycles}
One particular class of essentially left infinite vertices are those that {\em emit entries into cycles},
i.e. vertices $v$ that have some path $\alpha=e_1e_2\dots e_m$ of positive length, with $s(\alpha)=v$,
such that $e_1$ is an entry to a cycle.
Indeed, if $e_1$ enters a cycle $\nu$, then all paths $\nu^n\alpha$, $n\in \mathbb{N}$, are mutually 
incomparable.

Another class of essentially left infinite vertices are those that emit paths to infinitely many vertices. (In 
\cite{Tomforde4}, such vertices are called {\em left infinite}.)
\end{remark}

The following result generalizes \cite[Lemma 3.3(i)]{PaskRen1} and part of the proof of \cite[Theorem 3.2]{Tomforde4}. 

\begin{proposition}\label{tr-infinite}
Let $E$ be a directed graph, $g$ be a graph trace on $E$, and $v\in E^0$ be some vertex.
Assume either one of the hypotheses below is satisfied
\begin{itemize}
\item[$(a)$] $v$ emits an entry to a cycle; or
\item[$(b)$] $g$ is finite and $v$ is essentially left infinite.
\end{itemize}
Then $g(v)=0$.
\end{proposition}

\begin{proof}
The main ingredient in the proof is the observation that, for any finite set $F$ of mutually incomparable paths starting at $v$, one has the inequality
\begin{equation}
\sum_{w\in r(F)}g(w)\geq |F|\cdot g(v).
\label{tr-inf-obs}
\end{equation}
Indeed, if we list $F$ as $\{\alpha_1,\dots,\alpha_n\}$ (with all $\alpha$'s distinct, i.e. $n=|F|$), then by mutual incomparability, we have the inequality
$
\sum_{w\in r(F)}p^{}_w\geq\sum_{j=1}^np^{}_{\alpha_j}$, and then \eqref{tr-inf-obs} follows immediately from 
Theorem~\ref{gr-tr-char}.

By assumption, in either case, we can find an infinite set $Y\subset E^*$ of mutually incomparable paths starting at $v$,
such that the sum $M=\sum_{w\in r(Y)}g(w)$ is finite. (In case $(a)$, as seen in the preceding remark, we can ensure that $r(Y)$ is a singleton; case $(b)$ is trivial, by finiteness of $g$.)
The desired conclusion now follows immediately from \eqref{tr-inf-obs}, which implies $M\geq n\cdot g(v)$ for arbitrarily large $n$.
\end{proof}

\begin{mycomment}
As we will see shortly, graph traces on $E$ correspond to certain maps on the ``compactly supported'' diagonal subalgebra 
$\mathcal{D}(E)_{\text{fin}}=\bigcup_{V\in\mathcal{P}_{\text{fin}}(E^0)}\mathcal{D}(E)q^{}_V$, which will eventually yield tracial positive functionals on the dense $*$-subalgebra $C^*(E)_{\text{fin}}\subset C^*(E)$.
Although neither $\mathcal{D}(E)_{\text{fin}}$, nor $\mathcal{M}(E)_{\text{fin}}$, nor $C^*(E)_{\text{fin}}$,
 are $C^*$-algebras, they are nevertheless unions of increasing nets of unital $C^*$-algebras: 
$\mathcal{D}(E)_{\text{fin}}=\bigcup_{V\in\mathcal{P}_{\text{fin}}(E^0)}\mathcal{D}(E)q^{}_V$,
$\mathcal{M}(E)_{\text{fin}}=\bigcup_{V\in\mathcal{P}_{\text{fin}}(E^0)}\mathcal{M}(E)q^{}_V$,
and
$C^*(E)_{\text{fin}}=\bigcup_{V\in\mathcal{P}_{\text{fin}}(E^0)}q^{}_V C^*(E)q^{}_V$.
(Recall that, for any finite subset $V\subset E^0$, the projection $q^{}_V$ is defined to be $\sum_{v\in V}p^{}_v$.)
It is clear that the conditional expectations $\mathbb{E}_{\mathcal{M}}$ and $\mathbb{E}_{\mathcal{D}}$ map
$C^*(E)_{\text{fin}}$ onto $\mathcal{M}(E)_{\text{fin}}$ and
$\mathcal{D}(E)_{\text{fin}}$, respectively, so Corollaries \ref{cor-M-inv} and \ref{cor-D-inv} have suitable statements
applicable to $C^*(E)_{\text{fin}}$, with the word ``state'' replaced by ``positive linear functional.''
By definition,
positivity for linear functionals defined on each one of these $*$-algebras is equivalent to the positivity of their restrictions to each of the cut-off algebras corresponding to $V\in\mathcal{P}_{\text{fin}}(E^0)$.
Upon identifying
$\mathcal{D}(E)_{\text{fin}}=C_c(\widehat{\mathcal{D}(E)})$ and
$\mathcal{M}(E)_{\text{fin}}=C_c(\widehat{\mathcal{M}(E)})$, the positive cones
$\mathcal{D}(E)^+_{\text{fin}}$ and $\mathcal{M}(E)^+_{\text{fin}}$ correspond precisely to the non-negative 
continuous compactly supported functions. 
\end{mycomment}

With this set-up in mind, Theorem \ref{gr-tr-char} has the following consequence.

\begin{theorem}\label{gr-tr-lin}
For any graph trace $g$ on $E$, there exists a unique positive linear functional
$\eta=\eta_g:\mathcal{D}(E)_{\text{\rm fin}}\to\mathbb{C}$, such that
\begin{equation}
\eta_g(p^{}_\lambda)=g(s(\lambda)), \,\,\,\forall\,\lambda\in E^*.
\label{etag=}
\end{equation}
When restricted to the unital $C^*$-algebras $\mathcal{D}(E)q^{}_V$,
$V\in\mathcal{P}_{\text{\rm fin}}(E^0)$, the positive linear functionals
$\eta_g$, $g\in T(E)$, have norms:
$$\left\|\eta_g|_{\mathcal{D}(E)q^{}_V}\right\|=\sum_{v\in V}g(v).$$
In particular, for $g\in T(E)$, the functional $\eta_g$ is norm-continuous, if and only if
$g$ is finite, and in this case, one has $\|\eta_g\|=\|g\|_1$.
\end{theorem}

\begin{proof}

Let $\mathcal{A}$ be the complex span of $\{p_\lambda \}_{\lambda \in E^*}$, and let $\mathcal{A}_h$ be its Hermitean part, which is the same as the real span of $\{p_\lambda\}_{\lambda \in E^*}$. An application of Theorem \ref{gr-tr-char} shows that there is a unique $\mathbb{R}$-linear functional $\theta: \mathcal{A}_h \to \mathbb{R}$ with $\theta(p_\lambda) = g(s(\lambda))$ for all $\lambda \in E^*$. If we fix $V \in \mathcal{P}_{\operatorname{fin}}(E^0)$ and $x \in \mathcal{A}_h q_V$, another application of Theorem \ref{gr-tr-char} \crytserb{to the inequality $-||x|| q_V \leq x \leq ||x|| q_V$} shows that $|\theta(x)| \leq \theta(q_V) ||x||$. Thus for each $V \in \mathcal{P}_{\operatorname{fin}}(E^0)$, there is a unique $\mathbb{C}$-linear hermitean functional $\eta_V: \mathcal{D}(E) q_V \to \mathbb{C}$ with $||\eta_V|| = \eta_V(q_V)$, so that $\eta_V$ is in fact positive with norm equal to $\sum_{v \in V} g(v)$. Clearly if $V \subset W$ are both finite subsets of $E^0$, then $\eta_W|_{\mathcal{D}(E)q_V} = \eta_V$; thus, by density, there exists a unique positive linear functional $\eta_g$ defined on all of $\mathcal{D}(E)$ such that $\eta_g|_{\mathcal{D}(E)q_V} = \eta_V$ if $V \in \mathcal{P}_{\operatorname{fin}}(E^0)$.

\end{proof}

\begin{mycomment}
As a $*$-subalgebra in $C^*(E)_{\text{fin}}$, both $\mathcal{D}(E)_{\text{fin}}$ and
$\mathcal{M}(E)_{\text{fin}}$ are non-degenerate (since they both contain
$\{q^{}_V\}_{V\in\mathcal{P}_{\text{fin}}(E)}$, as well as regular, because they are normalized by
all $s^{}_e$, $e\in E^1$ and all $p^{}_v$, $v\in E^0$. Given a positive linear functional
$\eta$ on either one of these algebras, it then makes sense to define what it means for it to be $s^{}_e$-invariant.
\end{mycomment}

\begin{remark}
The map $g\longmapsto \eta_g$ establishes a affine bijective correspondence between
$T(E)$ and the space of positive linear functionals on $\mathcal{D}(E)_{\text{\rm fin}}$ that are
$s^{}_e$-invariant for all $e\in E^1$. The inverse of this correspondence is obtained as follows.
Given a linear positive functional $\theta$ on $\mathcal{D}(E)_{\text{fin}}$ which is $s^{}_e$-invariant, 
for all $e\in E^1$, the associated graph trace is simply the map 
\begin{equation}
g^\theta:E^0\ni v\longmapsto \theta(p^{}_v)\in [0,\infty).
\label{g-theta-def}
\end{equation}
When we specialize to the case of interest to us, Theorem \ref{gr-tr-lin} yields the following statement.
\end{remark}

\begin{theorem}\label{gr-tr-thm0}
For any normalized graph trace $g$, there exists a unique state $\psi_g\in S(\mathcal{D}(E))$ satisfying
\begin{equation}
\psi_g(p^{}_\lambda)=g(s(\lambda)), \,\,\,\forall\,\lambda\in E^*.
\label{psig=}
\end{equation}
All states $\psi_g$, $g\in T_1(E)$ are fully invariant, and furthermore, the correspondence
\begin{equation}
T_1(E)\ni g\longmapsto \psi_g\in S^{\text{\rm inv}}(\mathcal{D}(E))
\label{tr-to-Sinv}
\end{equation}
is an affine bijection, which has as its inverse the correspondence 
\begin{equation}
S^{\text{\rm inv}}(\mathcal{D}(E))\ni \theta\longmapsto g^\theta\in T_1(E)
\label{Sinv-tr}
\end{equation}
defined as in \eqref{g-theta-def}.
\qed
\end{theorem}

\begin{mycomment}
Using Corollary \ref{cor-D-inv}, 
it follows that for any $g\in T_1(E)$, the composition
$\chi_g=\psi_g\circ \mathbb{E}_{\mathcal{D}}$ defines a tracial state on $C^*(E)$; this way
we obtain an injective correspondence
\begin{equation}
T_1(E)\ni g\longmapsto \chi_g\in T(C^*(E)).
\label{gtr-to-tr}
\end{equation}
Of course, any tracial state $\tau\in T(C^*(E))$ becomes invariant, when restricted to $\mathcal{D}(E)$, so
using \eqref{Sinv-tr} we obtain a correspondence
\begin{equation}
 T(C^*(E)) \ni \tau \longmapsto g^\tau \in T_1(E). \label{surj}
\end{equation} 
Theorem \ref{gr-tr-thm0} shows that this map is surjective, because the correspondence \eqref{gtr-to-tr} is clearly an affine right inverse for \eqref{surj}.
The surjectivity of \eqref{surj} is also proved in \cite{Tomforde1}, by completely different means. 
\end{mycomment}
 
\begin{remark}\label{gauge-inv-tr=}
Using formulas \eqref{pD}, given a normalized graph trace $g\in T_1(E)$, the associated tracial state
$\chi_g=\psi_g\circ\mathbb{E}_{\mathcal{D}}$ -- hereafter referred to as the {\em Haar trace induced by $g$} -- 
acts on the spanning monomials as:
\begin{equation}
\chi_g(s^{}_\alpha s^*_\beta)=
\begin{cases}
g(s(\alpha)),&\text{if $\alpha=\beta$}\\
0,&\text{otherwise}
\end{cases}
\label{chig=}
\end{equation}
Among other things, the above formulas prove that $\chi_g$ is in fact {\em gauge invariant}, i.e.
$\chi_g\circ\gamma_z=\chi_g$, for all $z\in\mathbb{T}$. 

Conversely, every gauge invariant tracial state $\tau\in T(C^*(E))$ arises this way. Indeed, if $\tau$ is such a trace, then by gauge invariance it follows that, whenever $\alpha,\beta\in E^*$ are such that $|\alpha|\neq|\beta|$, we must have $\tau(s^{}_\alpha s^*_\beta)=0$; furthermore, if $|\alpha|=|\beta|$, then 
$$\tau(s^{}_\alpha s^*_\beta)=\tau(s^*_\beta s^{}_\alpha)=
\begin{cases}
\tau(0)=0,&\text{if $\alpha\neq\beta$}\\
\tau(s^*_\alpha s^{}_\alpha)=\tau(p^{}_{s(\alpha)}),&\text{otherwise}
\end{cases}
$$
so in all cases we get $\tau(s^{}_\alpha s^*_\beta)=\chi_{g^\tau}(s^{}_\alpha s^*_\beta)$.

To summarize:
\begin{itemize}
\item the range of the injective correspondence \eqref{gtr-to-tr} is the set
$T(C^*(E))^{\mathbb{T}}$ of gauge invariant tracial states; 
\item when restricting the correpondence
\eqref{surj} to $T(C^*(E))^{\mathbb{T}}$, one obtains an affine {\em isomorphism}
\begin{equation}
 T(C^*(E))^{\mathbb{T}} \ni \tau \longmapsto g^\tau \in T_1(E). \label{surj-gauge}
\end{equation}
\end{itemize} 
\end{remark}

When searching for an analogue of Theorem \ref{gr-tr-thm0}, with $\mathcal{D}(E)$ replaced by
$\mathcal{M}(E)$, it is obvious that the space $T(E)$ is not sufficient, so
additional structure needs to be added to it.

\begin{definition} \label{tagging-defn}
The {\em cyclic support\/} of a function $g:E^0\to \mathbb{C}$ is defined to \crytserb{be} the set
$$\text{supp}^cg=\{v\in E^0\,:\text{ $v\in E^0$ cyclic, $g(v)\neq 0$}\}.$$
(Recall that a cyclic vertex $v$ is one visited by a simple entry-less cycle. Equivalently,
$v$ is a \crytserb{ray of length zero}.) 
A \emph{cyclically tagged graph trace} consists of a pair $(g,\mu)$, where $g$ is a graph trace and
a map 
$\mu:\text{supp}^cg\ni v \longmapsto \mu_v\in
\operatorname{Prob}(\mathbb{T})$ -- hereafter referred to as the {\em tag\/}.
Note that our definition includes the possibility of an {\em empty\/} tag in the case 
when $\text{supp}^cg=\varnothing$. 
(More on this in Theorem~\ref{auto-gauge} below.)
The space of all such pairs will be denoted by $T^{\text{\sc ct}}(E)$.
The adjective ``finite,'' ``infinite,'' or ``normalized,'' is attached to $(g,\mu)$ precisely when it applies to $g$. 
\end{definition}

Using this terminology, one has the following extension of Theorem \ref{gr-tr-lin}.

\begin{theorem}\label{taged-gr-tr-lin}
For any cyclically tagged graph trace $(g,\mu)$ on $E$, there exists a unique positive linear functional
$\tilde{\eta}=\tilde{\eta}_{(g,\mu)}:\mathcal{M}(E)_{\text{\rm fin}}\to\mathbb{C}$, such that
\begin{itemize}
\item[(i)] $\tilde{\eta}_{(g,\mu)}(p^{}_\lambda)= g(s(\lambda))$, for every finite path $\lambda\in E^*$;
\item[(ii)] for any ray $\alpha$ and any integer $m\neq 0$,
$$\tilde{\eta}_{(g,\mu)}(b_\alpha ^m)=
\begin{cases}
g(s(\alpha))\int_{\mathbb{T}}z^m\,d\mu_{s(\alpha)}(z),&\text{ if }g(s(\alpha))\neq 0,\\
0,&\text{ otherwise}\end{cases}
$$
\end{itemize}
When restricted to the unital $C^*$-algebras $\mathcal{M}(E)q^{}_V$,
$V\in\mathcal{P}_{\text{\rm fin}}(E^0)$, the positive linear functionals
$\tilde{\eta}_{(g,\mu)}$, $(g,\mu)\in T^{\text{\rm\sc ct}}(E)$, have norms:
$$\left\|\tilde{\eta}_{(g,\mu)}|_{\mathcal{M}(E)q^{}_V}\right\|=\sum_{v\in V}g(v).$$
In particular, for any $(g,\mu)\in T^{\text{\rm\sc ct}}(E)$,
the functional $\tilde{\eta}_{(g,\mu)}$ is norm-bounded if and only if
$g$ is finite, and in this case, one has $\|\eta_{(g,\mu)}\|=\|g\|_1$.
\end{theorem}

\begin{proof}
Assume $(g,\mu)\in T^{\text{\rm\sc ct}}(E)$ is fixed throughout the entire proof.
Fix for the moment some a ray $\alpha$ with $g(s(\alpha))\neq 0$, and consider the $C^*$-subalgebra
$C^*(b_\alpha)\subset\mathcal{M}(E)$. (Recall that, if $\nu$ is the seed of the ray $\alpha$, then $b^{}_\alpha$ is the normal partial isometry $s^{}_\alpha s^{}_\nu s^*_\alpha$.) As pointed out in Lemma \ref{top-spec}, using the fact that the projection
$b_\alpha^0=p^{}_\alpha$ is the characteristic function of the compact-open set $\mathbf{T}_\alpha\subset\widehat{\mathcal{M}(E)}$, we have of course the equality
$\mathcal{M}(E)p^{}_\alpha=C^*(b^{}_\alpha)$,
so using the surjective $*$-homomorphism
$$\pi_\alpha:\mathcal{M}(E)\ni a\longmapsto ap^{}_\alpha \in C^*(b^{}_\alpha)\xrightarrow{\,\,\sim\,\,} 
C(\mathbb{T}),$$
we can define a state $\omega_\alpha$ on $\mathcal{M}(E)$ by
$$\omega_\alpha(a)=\int_{\mathbb{T}}\pi_\alpha(a)\,d\mu_{s(\alpha)}.$$
Specifically, if we write the compression $ap^{}_\alpha $ as a $f(b^{}_\alpha)$, for some 
$f\in C(\mathbb{T})$, then $\omega_\alpha(a)=\int_{\mathbb{T}}f(z)\,d\mu_{s(\alpha)}(z)$.
Using the product rules \eqref{GD-prod}, \eqref{GM-prod1} and \eqref{GM-prod2}, it follows that on the generator set
$G_{\mathcal{M}}(E)$, the state $\omega_\alpha$ acts
as
\begin{equation}
\omega_\alpha(p^{}_{\lambda})=
\begin{cases}
1,&\text{if }\lambda\prec\xi_\alpha;\\
0,&\text{otherwise;}\end{cases}
\quad
\omega_\alpha(b^m_{\alpha_1})=
\begin{cases}
\int_{\mathbb{T}}z^m\,d\mu_{s(\alpha)}(z),&\text{if }\alpha_1=\alpha;\\
0,&\text{otherwise.}
\end{cases}
\label{omega-alpha=}
\end{equation}
Define now the functional $\theta:\mathcal{M}(E)_{\text{fin}}\to\mathbb{C}$ by
\begin{equation}
\theta(a)=\sum_{\substack{\alpha\in E^*_{\text{\rm\sc ip}}\\ g(s(\alpha))\neq 0}}
g(s(\alpha))\omega_\alpha(a),\,\,\,a\in\mathcal{M}(E)_{\text{fin}}.
\label{theta=}
\end{equation}
Concerning the point-wise convergence of the sum in \eqref{theta=}, as well as its positivity, they are a consequence of the following fact.

\begin{claim}
For any vertex $v\in E^0$, one has the inequality
\begin{equation}
\sum_{\substack{\alpha\in E^*_{\text{\rm\sc ip}}\\ r(\alpha)=v}}g(s(\alpha))\leq g(v).
\label{theta-summable-claim}
\end{equation}
In particular, the sum
\begin{equation}
\theta_v=\sum_{\substack{\alpha\in E^*_{\text{\rm\sc ip}}\\ r(\alpha)=v}}g(s(\alpha))\omega_\alpha|_{\mathcal{M}(E)p^{}_v}
\label{theta-summable-claim0}
\end{equation}
is a norm-convergent sum, thus $\theta_v$ is a positive linear functional on $\mathcal{M}(E)p^{}_v$ with norm
\begin{equation}
\left\|\theta_v\right\|=\sum_{\substack{\alpha\in E^*_{\text{\rm\sc ip}}\\ r(\alpha)=v}}g(s(\alpha)).
\label{theta-summable-claim1}
\end{equation}
\end{claim}
The inequality \eqref{theta-summable-claim} follows from the observation that, for any finite set $F$ of rays with range $v$, the projections $\{p^{}_\alpha\}_{\alpha\in F}$ satisfy the inequality $\sum_{\alpha\in F}p^{}_\alpha\leq p^{}_v$, which by Theorem~\ref{gr-tr-char} implies
$\sum_{\alpha\in F}g(s(\alpha))\leq g(v)$. The equality \eqref{theta-summable-claim1} is now clear from the positivity
of $\theta_v$, which combined with \eqref{omega-alpha=} yields: 
$$\left\|\theta_v\right\|=\theta_v(p^{}_v)=
\sum_{\substack{\alpha\in E^*_{\text{\rm\sc ip}}\\ r(\alpha)=v}}g(s(\alpha))\omega_\alpha(p^{}_v)=
\sum_{\substack{\alpha\in E^*_{\text{\rm\sc ip}}\\ r(\alpha)=v}}g(s(\alpha)).$$

Using the Claim, we see that $\theta$ given in \eqref{theta=} is indeed correctly defined, positive and it can alternatively be presented as $\theta(a)=\sum_{v\in E^0}\theta_v(a)$ (a sum which has only finitely many non-zero terms for each $a\in\mathcal{M}(E)_{\text{fin}}$). By construction, $\theta$ acts on the generator set $G_{\mathcal{M}}(E)$ as:
\begin{align}
\theta(p^{}_\lambda)&=
\sum_{\substack{\alpha\in E^*_{\text{\rm\sc ip}}\\ \lambda\prec \xi_\alpha}}g(s(\alpha)),\,\,\,\lambda\in E^*
\label{theta-p}\\
\theta(b^m_{\alpha})&=
\begin{cases}
g(s(\alpha))\int_{\mathbb{T}}z^m\,d\mu_{s(\alpha)}(z),&\text{if }\alpha\in E^*_{\text{\rm\sc ip}}\text{ and }g(s(\alpha))\neq 0\\
0,&\text{otherwise}
\end{cases}
\label{theta-b}
\end{align}
Next we consider the positive linear functional $\eta_g:\mathcal{D}(E)_{\text{fin}}\to\mathbb{C}$ associated to $g$, as constructed in Theorem~\ref{gr-tr-lin}, and the linear positive functional 
$\eta_g\circ\mathbb{E}_{\mathcal{D}}:\mathcal{M}(E)_{\text{fin}}\to\mathbb{C}$.
(Here we use the fact that $\mathbb{E}_{\mathcal{D}}$ maps $C^*(E)_{\text{fin}}$ onto $\mathcal{D}(E)_{\text{fin}}$.)
Using Riesz' Theorem, there is a positive Radon measure $\upsilon$ on $\widehat{\mathcal{M}(E)}$, such that
$
\eta_g\left(\mathbb{E}_{\mathcal{D}}(f)\right)=
\int_{\widehat{\mathcal{M}(E)}}f\,d\upsilon$, for all $f\in C_c(\widehat{\mathcal{M}(E)})=
\mathcal{M}(E)_{\text{fin}}$. Using this measure, we now define the desired 
positive linear functional 
$\tilde{\eta}$ on $C_c(\widehat{\mathcal{M}(E)})=
\mathcal{M}(E)_{\text{fin}}$ by:
\begin{align}
\tilde{\eta}(f)&=
\theta(f)+\int_{\widehat{\mathcal{M}(E)}\smallsetminus\Omega_{\text{\rm\sc ip}}}f\,d\upsilon
=\notag\\
&=\theta(f)+\eta_g\left(\mathbb{E}_{\mathcal{D}}(f)\right)-
\sum_{\alpha\in E^*_{\text{\rm\sc ip}}}\eta_g\left(\mathbb{E}_{\mathcal{D}}(fp^{}_\alpha)\right)=
\label{eta-g-m}\\
&=\theta(f)+\eta_g\left(\mathbb{E}_{\mathcal{D}}(f)\right)-
\sum_{\alpha\in E^*_{\text{\rm\sc ip}}}\eta_g\left(\mathbb{E}_{\mathcal{D}}(f)p^{}_\alpha\right).
\label{eta-g-m0}
\end{align}
(The equality \eqref{eta-g-m} follows from Lemma \ref{top-spec}.)

To check condition (i), start with some $\lambda\in E^*$ and observe that, for all rays $\alpha$, 
we have the equalities
$$p^{}_\lambda p^{}_\alpha=\begin{cases}
p^{}_\alpha,&\text{if }\lambda\prec \xi_\alpha\\
0,&\text{otherwise}
\end{cases}
$$
which by \eqref{theta-p} imply that
$$
\sum_{\alpha\in E^*_{\text{\rm\sc ip}}}\eta_g\left(\mathbb{E}_{\mathcal{D}}(p^{}_\lambda p^{}_\alpha)\right)
=
\sum_{\substack{\alpha\in E^*_{\text{\rm\sc ip}}\\ \lambda\prec\xi_\alpha}}\eta_g(p^{}_\alpha)=
\sum_{\substack{\alpha\in E^*_{\text{\rm\sc ip}}\\ \lambda\prec\xi_\alpha}}g(s(\alpha))=\theta(p^{}_\lambda),$$
so by \eqref{eta-g-m} we obtain the desired property
$$\tilde{\eta}(p^{}_\lambda)=\eta_g(p^{}_\lambda)=g(s(\lambda)).$$

In order to check condition (ii), we simply verify that, for any ray $\alpha$ and any integer $m$, we have the equality \begin{equation}
\tilde{\eta}(b^m_\alpha)=\theta(b^m_\alpha).
\label{cond(ii)}
\end{equation} The case when $m=0$ we have $b^0_\alpha=p^{}_\alpha$,
so by condition (i) and
\eqref{theta-b}, we have $\tilde{\eta}(b^0_\alpha)=\tilde{\eta}(p^{}_\alpha)=g(s(\alpha))=\theta(b^0_\alpha)$. In the case when $m\neq 0$, we notice that since $\mathbb{E}_{\mathcal{D}}$ vanishes on 
$G(E)\smallsetminus G_{\mathcal{D}}(E)$ -- by \eqref{pD} -- we have
$\mathbb{E}_{\mathcal{D}}(b^m_\alpha)=0$, and then \eqref{cond(ii)} is trivial using
\eqref{eta-g-m0}.

The remaining statements in the Theorem (including the uniqueness of $\tilde{\eta}$) are pretty clear, since any
positive linear functional $\tilde{\eta}$ satisfying conditions (i) and (ii) must satisfy
$\tilde{\eta}|_{\mathcal{D}(E)_{\text{fin}}}=
\eta_g$, from which the continuity of the restrictions 
$\tilde{\eta}|_{\mathcal{M}(E)q^{}_V}$ follows immediately.
\end{proof}

One aspect not addressed so far is invariance of the states $\tilde{\eta}$. For this purpose, the following definition is well-suited.

\begin{definition}
Two cyclic vertices are said to be \emph{equivalent} if they are visited by the same entry-less cycle. A cyclically tagged graph trace $(g,\mu)$ is said to be {\em consistent} if $\mu_v = \mu_{v'}$ whenever $v$ and $w$ are equivalent.  
(Note that if two cyclic vertices $v,w$ are equivalent, then $g(v)=g(w)$.)
The space of all consistent cyclically tagged traces on $E$ is denoted by $T^{\text{\rm\sc cct}}(E)$.
As agreed earlier, the adjective ``finite,'' ``infinite,'' or ``normalized,'' is attached to an element
$(g,\mu)\in T^{\text{\rm\sc cct}}(E)$, precisely when it applies to $g$. In particular, the space of normalized consistent cyclically tagged graph traces on $E$ is denoted by
$T_1^{\text{\rm\sc cct}}(E)$.
\end{definition}

\begin{proposition}\label{cct-equiv}
A cyclically tagged graph trace $(g,\mu)$ is consistent if and only if the associated positive functional $\tilde{\eta}_{(g,\mu)}:\mathcal{M}(E)_{\text{fin}} \to \mathbb{C}$ constructed in Theorem~\ref{taged-gr-tr-lin} is $s^{}_e$-invariant for all $e \in E^1$.
\end{proposition}
\begin{proof}
Assume $(g,\mu)$ is consistent, and let us show the invariance of $\tilde{\eta}_{(g,\mu)}$,
which amounts to checking, that for each  $e\in E^1$, we have:
\begin{itemize}
\item[(i)] 
$\tilde{\eta}_{(g,\mu)}(s^{}_ep^{}_\lambda s^*_e)=
\tilde{\eta}_{(g,\mu)}(p^{}_ep^{}_\lambda)$, $\forall\,\lambda\in E^*$;
\item[(ii)] 
$\tilde{\eta}_{(g,\mu)}(s^{}_eb^m_\alpha s^*_e)=
\tilde{\eta}_{(g,\mu)}(p^{}_eb^m_\alpha)$, $\forall\,\alpha\in E^*_{\text{\sc ip}}$, $m\in\mathbb{Z}$.
\end{itemize}
Property (i) is obvious, since $\tilde{\eta}_{(g,\mu)}$ agrees with the $s^{}_e$-invariant functional $\eta_g$ on $\mathcal{D}(E)_{\text{fin}}$. As for condition (ii), we only need to verify it if 
$s(e)=r(\alpha)$ (otherwise both sides are zero). Also notice that if $|\alpha|>0$, then $e\alpha$ is also a ray 
with $s(e\alpha)=s(\alpha)$, which satisfies $s^{}_eb^m_\alpha s^*_e=b^m_{e\alpha}$, so by condition (ii) in Theorem~\ref{taged-gr-tr-lin}, we have $\tilde{\eta}_{(g,\mu)}(s^{}_eb^m_\alpha s^*_e)=
\tilde{\eta}_{(g,\mu)}(b^m_{e\alpha})=g(s(e\alpha))\int_{\mathbb{T}}z^m\,d\mu_{s(e\alpha)}(z)=
g(s(\alpha))\int_{\mathbb{T}}z^m\,d\mu_{s(\alpha)}(z)=
\tilde{\eta}_{(g,\mu)}(b^m_\alpha)$.
In the remaining case, $|\alpha|=0$, so $\alpha$ reduces to a vertex $v=r(\nu)$, for some simple entry-less cycle $\nu$. If $e$ is not an edge in $\nu$, then it is a ray, thus the preceding argument still applies 
(we will have $s^{}_eb^m_vs^*_e=b^m_e$).
If $e$ is an edge on $\nu$, then $s^{}_eb^m_vs^*_e=b^m_{r(e)}$, with $r(e)$ obviously equivalent to $v$, and the desired equality -- which now reads
$\tilde{\eta}_{(g,\mu)}(b^m_{r(e)})=
\tilde{\eta}_{(g,\mu)}(b^m_v)$ -- follows from the equalities $g(v)=g(r(e))$ and  $\mu_v=\mu_{r(e)}$.

Conversely, notice first that, if $\tilde{\eta}_{(g,\mu)}$ is $s^{}_e$-invariant, for all $e\in E^1$, then it will also satisfy
the identity
\begin{equation}
\tilde{\eta}_{(g,\mu)}(s^{}_\lambda a s^*_\lambda)=
\tilde{\eta}_{(g,\mu)}(p^{}_\lambda a),\,\,\,\forall\,\lambda\in E^*, a\in \mathcal{M}(E)_{\text{fin}}.
\label{tilde-eta-ainv}
\end{equation}
Secondly, observe that, if $v$, $v'$ are equivalent cyclic vertices, presented as $v=s(\nu)$ and
$v'=s(\nu')$ for two simple entry-less cycles, then we can write $\nu=\alpha\beta$ and $\nu'=\beta\alpha$ for two suitably chosen paths $\alpha,\beta\in E^*$. This clearly implies that
$b^{}_{v'}=s^{}_\beta b^{}_v s^*_\beta$, which also yields
$b^m_{v'}=s^{}_\beta b^m_v s^*_\beta$, $\forall\,m\in\mathbb{Z}$.

Combining these two observations with condition (ii) from Theorem~\ref{taged-gr-tr-lin}, it follows that, if $\tilde{\eta}_{(g,\mu)}$ is invariant, then for any two equivalent cyclic vertices $v$ and $v'$ we have (with $\alpha$, $\beta$ as above):
\begin{align*}
\int_{\mathbb{T}}z^m\,d\mu_{v'}(z)&=
\tilde{\eta}_{(g,\mu)}(b^m_{v'})=
\tilde{\eta}_{(g,\mu)}(s^{}_\beta b^m_{v}s^*_\beta)=
\tilde{\eta}_{(g,\mu)}(p^{}_\beta b^m_{v})=\\
&=
\tilde{\eta}_{(g,\mu)}(b^m_{v})=
\int_{\mathbb{T}}z^m\,d\mu_{v}(z),
\,\,\,\forall\,m\in\mathbb{Z},
\end{align*}
which clearly implies $\mu_{v'}=\mu_v$.
\end{proof}

\begin{remark}
The map $(g,\mu)\longmapsto \tilde{\eta}_{(g,\mu)}$ establishes a affine bijective correspondence between
$T^{\text{\sc cct}}(E)$ and the space of positive linear functionals on $\mathcal{M}(E)_{\text{\rm fin}}$ that are
$s^{}_e$-invariant for all $e\in E^1$. The inverse of this correspondence is the map
$\theta\longmapsto (g^\theta,\mu^\theta)$ defined as follows.
Given a linear positive functional $\theta$ on $\mathcal{M}(E)_{\text{fin}}$ which is $s^{}_e$-invariant, 
for all $e\in E^1$, the graph trace $g^\theta$ is given by \eqref{g-theta-def}, and the tag
$\mu^\theta=(\mu^\theta_v)_{v\in\text{supp}^cg^\theta}$ is given (implicitly) by
\begin{equation}
\int_{\mathbb{T}}f(z)\,d\mu^\theta_v(z)=
\frac{\theta(f(b^{}_v))}{g^\theta(v)},\,\,\,\forall\,v\in\text{supp}^cg^\theta,\,f\in C(\mathbb{T}).
\label{im-def-muv}
\end{equation}
\end{remark}

When we specialize to states, we now have the following extension of Theorem~\ref{gr-tr-thm0}.

\begin{theorem}\label{gr-tr-thm1}
For any normalized consistent cyclically tagged graph trace $(g,\mu)\in T^{\text{\rm\sc cct}}_1(E)$, there exists a unique state $\phi_{(g,\mu)}\in S(\mathcal{M}(E))$ satisfying
\begin{itemize}
\item[(i)] $\phi_{(g,\mu)}(p^{}_\lambda)= g(s(\lambda))$, for every finite path $\lambda\in E^*$;
\item[(ii)] for any ray $\alpha$ and any integer $m$:
$$\phi_{(g,\mu)}(b_\alpha ^m)=
\begin{cases}
g(s(\alpha))\int_{\mathbb{T}}z^m\,d\mu_{s(\alpha)}(z),&\text{ if }g(s(\alpha))\neq 0,\\
0,&\text{ otherwise}\end{cases}
$$
\end{itemize}
All states $\phi_{(g,\mu)}$, $(g,\mu)\in T^{\text{\rm\sc cct}}_1(E)$ are fully invariant, and furthermore, the correspondence
\begin{equation}
T^{\text{\rm\sc cct}}_1(E)\ni (g,\mu)\longmapsto \phi_{(g,\mu)}\in S^{\text{\rm inv}}(\mathcal{M}(E))
\label{tr-to-SinvM}
\end{equation}
is an affine bijection, which has as its inverse the correspondence 
\begin{equation}
S^{\text{\rm inv}}(\mathcal{M}(E))\ni \theta\longmapsto (g^\theta,\mu^\theta)\in T^{\text{\rm\sc cct}}_1(E)
\label{Sinv-trM}
\end{equation}
defined as in \eqref{g-theta-def} and \eqref{im-def-muv}.
\qed
\end{theorem}

\begin{mycomment}
Using Corollary \ref{cor-M-inv}, it follows that for any $(g,\mu)\in T^{\text{\rm\sc cct}}_1(E)$, the composition
$\tau_{(g,\mu)}=\phi_{(g,\mu)}\circ \mathbb{E}_{\mathcal{M}}$ defines a tracial state on $C^*(E)$; this way
we obtain an injective correspondence
\begin{equation}
T^{\text{\rm\sc cct}}_1(E)\ni (g,\mu)\longmapsto \tau_{(g,\mu)}\in T(C^*(E)).
\label{cctgtr-to-tr}
\end{equation}
Of course, any tracial state $\tau\in T(C^*(E))$ becomes invariant, when restricted to $\mathcal{M}(E)$, so
using \eqref{Sinv-trM} we obtain a correspondence
\begin{equation}
 T(C^*(E)) \ni \tau \longmapsto (g^\tau,\mu^\tau) \in T^{\text{\rm\sc cct}}_1(E). \label{surjM}
\end{equation} 
Theorem \ref{gr-tr-thm1} shows that this map is surjective, because the correspondence \eqref{cctgtr-to-tr} is clearly an affine right inverse for \eqref{surjM}.
\end{mycomment}

\begin{remark}
The range of \eqref{cctgtr-to-tr} clearly contains the range of \eqref{gtr-to-tr}, which equals
$T(C^*(E))^{\mathbb{T}}$. After all, any trace $g\in T_1(E)$ can be tagged using the constant
map $\mu:\text{supp}^cg\to\text{Prob}(\mathbb{T})$ that takes $\mu_v$ to be the Haar measure for every $v$, and it is straightforward to verify that for this particular tagging one, has $\tau_{(g,\mu)}=\chi_g$.
\end{remark}

Concerning the range of \eqref{cctgtr-to-tr}, one legitimate question is whether it equals the whole tracial state space $T(C^*(E))$. Using the bijection \eqref{tr-to-SinvM}, this question is equivalent to the surjectivy of the map
\begin{equation}
S^{\text{inv}}(\mathcal{M}(E))\ni \phi\longmapsto \phi\circ\mathbb{E}_{\mathcal{M}}\in T(C^*(E)).
\label{SinvM-to-tr}
\end{equation}
As we have seen in Corollary~\ref{cor-Sinv-T-iso}, a sufficient condition for the surjectivity
of \eqref{SinvM-to-tr} is the condition that the inclusion $\mathcal{M}(E)\subset C^*(E)$ has the (honest) extension property. As it turns out, this issue can be neatly described using the graph.

\begin{theorem}\label{tight-thm}
The inclusion $\mathcal{M}(E) \subset C^*(E)$ has the extension property, if and only if
no cycle in $E$ has an entry. 
\end{theorem}

\begin{proof}
To prove the ``if'' implication, assume that no cycle in $E$ has an entry, fix a pure state $\omega$ on $\mathcal{M}(E)$, and let $\phi$ be an extension of $\omega$ to $C^*(E)$. In order to prove uniqueness of $\phi$, it suffices to show that the value of $\phi$ on a standard generator $s_\alpha^{} s_\beta^*$ is independent of the choice of $\phi$. By assumption, there is a $x \in E^{\leq \infty}$ and $z \in \mathbb{T}$ such that $\omega = \omega_{z,x}$ as in Lemma \ref{top-spec}. On the one hand, by Fact 3.1 and the observation that $\omega(p_\gamma)=1$ for all $\gamma \prec x$, it follows that 
\begin{equation} 
\label{tight-state} \forall \,\gamma \prec x:\quad\phi(s_\alpha^{} s_\beta^*)=\phi(p_\gamma^{} s_\alpha^{} s_\beta^* p_\gamma^{}).
\end{equation}
 On the other hand, using the results from \cite[Section 3]{NR}, it follows that there is $\gamma \prec x$ such that $p_\gamma^{} s_\alpha^{} s_\beta^* p_\gamma^{}$ belongs to $\mathcal{M}(E)$. (In the language of \cite{NR}, $x$ must be essentially aperiodic by our assumption on $E$.) Using \eqref{tight-state} it follows that $\phi(s_\alpha^{} s_\beta^*) = \omega(p_\gamma^{} s_\alpha^{} s_\beta^* p_\gamma^{})$, and the desired conclusion follows.

For the ``only if'' direction, we show that if there is a cycle $\nu\in E^*$ that has an entry, then we can construct a pure state on $\mathcal{M}(E)$ which has multiple extensions to states on $C^*(E)$. Consider the path $x=\nu^\infty\in E^\infty$ formed by following $\nu$ infinitely many times. For each $z\in\mathbb{T}$ consider the state
$\omega_{z,x}\in S(C^*(E))$ introduced in Definition~\ref{twisted-path-rep}, given by
$$\omega_{z,x}(a) = \langle \delta_x | \pi_{\text{path}}(\gamma_z(a)) \delta_x \rangle.$$
As explained in Remark~\ref{core-embed}, since $x\not\in E^\infty_{\text{\sc ip}}$, it follows that:
$$(z,x)\sim (1,x),\,\,\,\forall\,z\in\mathbb{T},$$
which by Lemma~\ref{top-spec} means that all restrictions $\omega_{z,x}|_{\mathcal{M}(E)}$, $z\in\mathbb{T}$, coincide, so they are all equal to the pure state $\vartheta\in\widehat{\mathcal{M}(E)}$ corresponding to the equivalence class 
$(1,x)_{\sim}=\mathbb{T}\times\{x\}$. However, as states on $C^*(E)$, the functionals $\omega_{z,x}$, $z\in\mathbb{T}$ cannot all be equal, since for example we have $\omega_{z,x}(\nu)=z^{|\nu|}$, $\forall\,z\in\mathbb{T}$.
\end{proof}

\begin{definition}
A graph $E$ is {\em tight}, if every cycle is entry-less.
\end{definition}

Combining Theorem~\ref{tight-thm} with Corollary~\ref{cor-Sinv-T-iso} and
Theorem~\ref{gr-tr-thm1} we now obtain the following statement.

\begin{theorem}\label{traces-on-tight}
If $E$ is tight, then the correspondence \eqref{cctgtr-to-tr} is an affine isomorphism between
the space $T^{\text{\sc cct}}_1(E)$ and the tracial state space $T(C^*(E))$.\qed
\end{theorem}
 
\begin{remark}
Tight graphs are interesting in other respects: they are the only graphs that yield finite, stably finite, quasi-diagonal, or AF-embeddable $C^*$-algebras (\cite{Schaf}), as well as the only graphs that yield graph algebras with stable rank one (\cite{JPS}). A graph which yields a $C^*$-algebra with Hausdorff spectrum must be tight, although this is not sufficient \cite[Ex. 10]{Goehle2}.
\end{remark}

In the remainder of this paper we aim to parametrize the entire tracial state space $T(C^*(E))$ for arbitrary graphs by employing Theorem~\ref{traces-on-tight} in conjunction with certain procedures
that replace the graph $E$ with a tight sub-graph $E'$, in such a way that
the tracial state spaces $T(C^*(E))$ and $T(C^*(E'))$ coincide. Since the sub-graphs that are best suited for analyzing 
how the trace spaces change are the {\em canonical\/} ones, the following terminology is all we need.


\begin{definition} \label{tightening-defn}
If $E$ is a directed graph, a {\em tightening\/} of $E$ is a canonical sub-graph, i.e. one that can be presented as $E\setminus H$, for some saturated hereditary subset $H\subset E^0$, in such a way that
\begin{itemize}
\item[{\sc (a)}] $E\setminus H$ is tight, and
\item[{\sc (b)}] the canonical $*$-homomorphism $\rho_H:C^*(E)\to C^*(E\setminus H)$ implements a bijective correspondence:
$T(C^*(E\setminus H))\ni \tau\longmapsto \tau\circ \rho_H\in T(C^*(E))$
\end{itemize}
Since $\rho_H$ is always surjective, the correspondence from {\sc (b)} is always injective, so the only requirement in our definition is its {\em surjectivity}.
\end{definition}

When it comes to parametrizing tracial states on graph $C^*$-algebras, the most useful and natural
tightening is as follows.

\begin{example}\label{min-tight}
Let $E$ be a graph, and let $C=C_E$ be the set of vertices which emit entrances into cycles. 
The set $C$ is obviously hereditary, but not saturated in general, so we need to take
its saturation $\overline{C}$.
As it turns out, $E\setminus \overline{C}$ constitutes a tightening of $E$. First of all, since passing from 
$E$ to $E\setminus \overline{C}$ clearly removes all entries into the cycles in $E$, it is clear that
$E\setminus \overline{C}$ is tight. Secondly, in order to justify the surjectivity of
\begin{equation}
T(C^*(E\setminus \overline{C}))\ni \tau\longmapsto \tau\circ \rho_{\overline{H}}\in T(C^*(E)),
\label{rho-traces-min}
\end{equation}
all we must show is the fact that {\em all tracial states on $C^*(E)$ vanish on}
$\ker\rho_{\overline{C}}$, for which it suffices to prove the inclusion
$H\subset N_g$, which in itself is a consequence of Proposition~\ref{tr-infinite}.
\end{example}

The sub-graph constructed in the above Example is called the {\em minimal tightening}, and is denoted by
$E_{\text{tight}}$. The canonical $*$-homomorphism will be denoted by $\rho_{\text{tight}}:C^*(E)\to C^*(E_{\text{tight}})$.
Combining this construction with Theorem~\ref{traces-on-tight} we now obtain.

\begin{theorem}
For any directed graph $E$, the map
$$T_1^{\text{\rm\sc cct}}(E_{\text{\rm tight}})\ni (g,\mu)\longmapsto 
\tau_{(g,\mu)}\circ \rho_{\text{\rm tight}}\in T(C^*(E))$$
is an affine isomorphism.\qed
\end{theorem} 

The final result in this paper deals with a graph-theoretic characterization of automatic gauge invariance for tracial states, which as pointed out in Remark~\ref{gauge-inv-tr=} is equivalent to the surjectivity of the map \eqref{gtr-to-tr}. 
In \cite{Tomforde4}, it is shown that this feature is implied by condition (K). However, as Theorem~\ref{auto-gauge}
below shown, this is not necessary. 

\begin{theorem}\label{auto-gauge}
For a directed graph $E$, the following conditions are equivalent:
\begin{itemize}
\item[(i)] all tracial states on $C^*(E)$ are gauge invariant;
\item[(ii)] the source of each cycle in $E$ is essentially left infinite.
\end{itemize}
\end{theorem}

\begin{proof}
(i) $\Rightarrow$ (ii): 
Suppose that  $\lambda = e_1 \ldots e_m$ is a cycle such that $v=s(\lambda)=r(e_1)$ is not essentially left infinite; we show how to construct a tracial state on $C^*(E)$ which is not gauge-invariant. Note that as $v$ is not essentially infinite, in particular it does not emit an entrance to any cycle; therefore, none of the edges in $\lambda$ will be removed when forming $E_{\operatorname{tight}}$, and so we can assume that $E$ is tight. (Since the canonical quotient $\pi: C^*(E) \to C^*(E_{\operatorname{tight}})$ is equivariant for the respective gauge actions, a non-gauge invariant tracial state on $C^*(E)_{\operatorname{tight}})$ will give rise to a non-gauge invariant trace on $C^*(E)$.) 

Say that a path $\mu \in E^*$ is \emph{acyclic} if it cannot be written as $\mu = \alpha \nu \beta$ for $\alpha,\beta \in E^*$ and $\nu$ a cycle. Let $A$ denote the set of all acyclic paths with source $v$; note that any two paths in $A$ are incomparable, and so $A$ must be finite because $v$ is not essentially left infinite. For $w \in E^0$ let $g(w) = |A \cap r^{-1}(w)|$; it is straightforward to verify that $g$ is a finite graph trace with $g(v) = 1$ which we can normalize to obtain $g' \in T_1(E)$. Note that the cyclic support of $g'$ is precisely $r(\{e_1,\ldots,e_m\})$ (as $v$ is not essentially left infinite, it emits no entrances to cycles).

Now we can take any $z \in \mathbb{T} \setminus \mathbb{U}_{| \lambda |}$ and let $\mu_{s(e_i)}=\delta_z$ for all $i=1,\ldots,m$. The affiliated tracial state $\tau_{(g,\mu)} \in T(C^*(E))$ will satisfy \[ \tau_{(g,\mu)}(b_\lambda) = g(s(\lambda)) z^{|\lambda|} \neq 0 \]so that in particular $\tau_{(g,\mu)}$ is not gauge-invariant. 

(ii) $\Rightarrow$ (i): Suppose that the source of each cycle is essentially left infinite. Any finite graph trace must vanish on an essentially left infinite vertex as in Proposition \ref{tr-infinite}; hence if every source of every cycle is essentially left infinite, then there are no vertices in the cyclic support of any graph trace, and so there are no taggings to consider. Thus every tracial state on $C^*(E_{\operatorname{tight}})$ is gauge-invariant, which  shows that every tracial state on $C^*(E)$ is gauge-invariant. 
\end{proof}

\begin{mycomment}
Besidese the minimal tightening $E_{\operatorname{tight}}$ introduced in this paper, other tightenings could naturally be considered. The same arguments as those used in Example~\ref{min-tight} can be used with $C$ replaced by another hereditary subset
$H\subset E^0$, as long as:
\begin{itemize}
\item[{\sc (a)}] the canonical sub-graph $E\setminus \overline{H}$ is tight, and
\item[{\sc (b)}] one has the inclusion $H\subset N_g$, for all $g\in T_1(E)$. 
\end{itemize}
One way to ensure {\sc (a)} is to take $H$ to contain $C_E$. As far as condition {\sc (b)} is concerned, we could use 
Proposition~\ref{tr-infinite} as a guide. In particular, we can consider the set
$L=L_E$ of \emph{all} essentially left infinite vertices. Since $L_E$ is potentially much larger than $C_E$, the resulting subgraph $E \setminus \overline{L}_E$ will potentially be considerably smaller than $E_{\operatorname{tight}}$ (and thus easier to analyze regarding graph traces). 
\end{mycomment}

\bibliographystyle{plain}
\bibliography{Bibliography2016-04-08}

\begin{thebibliography}{10}

\bibitem{ABG}
R.J. Archbold, J.W. Bunce, and K.D. Gregson.
\newblock Extensions of states of {$C^*$}-algebras, ii.
\newblock {\em Proc. Edinb. Math. Soc.}, 92A:113--122, 1982.

\bibitem{BNR}
J.~Brown, G.~Nagy, and S.~Reznikoff.
\newblock A generalized cuntz--krieger uniqueness theorem for higher-rank
  graphs.
\newblock {\em J. Funct. Anal.}, 266:2590--2609, 2013.

\bibitem{BNRSW}
J.~Brown, G.~Nagy, S.~Reznikoff, A.~Sims, and D.~Williams.
\newblock Cartan {S}ubalgebras in {$C^*$}-{A}lgebras of {H}ausdorff {E}tale
  {G}roupoids. arxiv:1503.03521.

\bibitem{Goehle2}
G.~Goehle.
\newblock Groupoid {$C^*$}-algebras with {H}ausdorff spectrum.
\newblock {\em Bull. Aust. Math. Soc.}, 88:232--242, 2013.

\bibitem{HuefRaeburn}
{A. an} Huef and I.~Raeburn.
\newblock The ideal structure of cuntz–krieger algebras.
\newblock {\em Ergodic Theory Dynam. Systems}, 17:611–624, 1997.

\bibitem{JPS}
J.~Jeong, G.~Park, and D.~Shin.
\newblock Stable rank and real rank of graph {$C^*$}-algebras.
\newblock {\em Pac. J. Math}, 200(2):331--343, 2001.

\bibitem{KS}
R.~Kadison and I.~Singer.
\newblock Extensions of pure states.
\newblock {\em Amer. Jour. Math}, 81(2):383--400, 1960.

\bibitem{Kumjian}
A.~Kumjian.
\newblock On {$C^*$}-{D}iagonals.
\newblock {\em Can. J. Math.}, 38(4):969--1008, 1986.

\bibitem{Lance}
{E. C.} Lance.
\newblock {\em {H}ilbert {$C^*$}-modules: {A} toolkit for operator
  algebraists}, volume 210 of {\em {L}ondon {M}ath {S}oc. {L}ecture {N}ote
  {S}eries}.
\newblock {C}ambridge {U}niv. {P}ress, 1994.

\bibitem{NR}
G.~Nagy and S.~Reznikoff.
\newblock Abelian core of graph algebras.
\newblock {\em J. London Math. Society}, 3:889--908, 2012.

\bibitem{NR2}
G.~Nagy and S.~Reznikoff.
\newblock Pseudo-diagonals and uniqueness theorems.
\newblock {\em Proc. Amer. Math. Soc.}, 142(1):263--275, January 2014.

\bibitem{PaskRen1}
D.~Pask and A.~Rennie.
\newblock The noncommutative geometry of graph ${C^*}$-algebras {I}: The index
  theorem.
\newblock {\em J. Funct. Anal.}, 233:92--134, 2006.

\bibitem{Raeburn}
I.~Raeburn.
\newblock {\em {G}raph {A}lgebras}.
\newblock CBMS Lecture Notes. American Mathematical Society, 2005.

\bibitem{Renault}
J.~Renault.
\newblock {\em A {G}roupoid {A}pproach to {$C^*$}-{A}lgebras}.
\newblock Lecture Notes in Mathematics. Springer Verlag, 1980.

\bibitem{Renault3}
J.~Renault.
\newblock Cartan {S}ubalgebras in {$C^*$}-{A}lgebras.
\newblock {\em Irish Math. Soc. Bulletin}, 61:29--63, 2008.

\bibitem{Schaf}
C.~Schafhauser.
\newblock Af-embeddings of {G}raph {A}lgebras. ar{X}iv:1405.7757.

\bibitem{Tomforde1}
M.~Tomforde.
\newblock The ordered ${K_0}$-group of a graph ${C^*}$-algebra.
\newblock {\em C.R. Math. Acad. Sci. Soc.}, 25:19--25, 2003.

\bibitem{Tomforde2}
M.~Tomforde.
\newblock A unified approach to {E}xel-{L}aca algebras and ${C^*}$-algebras
  associated to graphs.
\newblock {\em J. Operator Theory}, 50:345--368, 2003.

\bibitem{Tomforde4}
M.~Tomforde.
\newblock Stability of ${C^*}$-algebras affiliated to graphs.
\newblock {\em Proc. Amer. Math. Soc.}, 132:1787--1795, 2004.

\end{thebibliography}
\end{document}